\documentclass[%
  onecolumn
   , colorlinks 
]{mpi2015-cscpreprint}

\usepackage[american]{babel}

\usepackage{graphicx}

\usepackage{amssymb}
\usepackage{amsthm}

 \usepackage{multirow}
 \usepackage{booktabs}

\usepackage{subcaption}
\usepackage{graphicx}
\usepackage[capitalise]{cleveref}
\crefname{algocf}{alg.}{algs.}
\Crefname{algocf}{Algorithm}{Algorithms}

\AtBeginDocument{
  \newtheorem{theorem}{Theorem}[section]
  \newtheorem{lemma}[theorem]{Lemma}
  \newtheorem{definition}{Definition}[section]
  \newtheorem{remark}{Remark}[section]
  \newtheorem{property}{Property}[section]
}
\crefname{definition}{Definition}{Definitions}
\crefname{remark}{Remark}{Remarks}
\crefname{property}{Property}{Properties}

\usepackage[vlined,longend,linesnumbered,ruled]{algorithm2e}
\usepackage[round]{natbib}

\hypersetup{citecolor=mpigreen}

\usepackage[scr=boondox]{mathalfa}

\input{mymymacros.sty}

\begin{document}


\title{Discrete empirical interpolation in the\\ tensor t-product framework}
  
\author[$\ast$]{Sridhar Chellappa}
\affil[$\ast$]{Max Planck Institute for Dynamics of Complex Technical Systems, Magdeburg, Germany.\authorcr%
 \vspace{1mm}
  \email{chellappa@mpi-magdeburg.mpg.de}, \orcid{0000-0002-7288-3880} \authorcr
  \email{feng@mpi-magdeburg.mpg.de}, \orcid{0000-0002-1885-3269} \authorcr
  \email{benner@mpi-magdeburg.mpg.de}, \orcid{0000-0003-3362-4103} 
}
\author[$\ast$]{Lihong Feng}  
\author[$\ast$]{Peter Benner}

  
\shorttitle{t-tensor-product discrete empirical interpolation}
\shortauthor{}
\shortdate{}
  
\keywords{Function approximation, Discrete empirical interpolation, Sparse sampling, State reconstruction, Sensor selection}

  
\abstract{%
	The discrete empirical interpolation method (DEIM) is a well-established approach, widely used for state reconstruction using sparse sensor/measurement data, nonlinear model reduction, and interpretable feature selection. We introduce the tensor t-product Q-DEIM (t-Q-DEIM), an extension of the DEIM framework for dealing with tensor-valued data. The proposed approach seeks to overcome one of the key drawbacks of DEIM, viz., the need for matricizing the data, which can distort any structural and/or geometric information. Our method leverages the recently developed tensor t-product algebra to avoid reshaping the data. In analogy with the standard DEIM, we formulate and solve a tensor-valued least-squares problem, whose solution is achieved through an interpolatory projection. We develop a rigorous, computable upper bound for the error resulting from the t-Q-DEIM approximation. Using five different tensor-valued datasets, we numerically illustrate the better approximation properties of t-Q-DEIM and the significant computational cost reduction it offers.
  }

 \novelty{\begin{itemize}
 		\item A new extension of the discrete empirical interpolation method, called t-Q-DEIM, to deal with three dimensional tensor-structured data, is proposed
 		\item A rigorous upper bound for the error due to the t-Q-DEIM approximation of the tensor-valued data is developed
 		\item A numerical recipe for computing the necessary elements involved in this extended DEIM formulation is provided and it is tested on relevant examples
 \end{itemize}}

\maketitle

  
\section{Introduction}%
\label{sec:intro}
Approximating parametrized functions or matrices from sparsely sampled data is an important problem with applications in several fields. Within nonlinear model order reduction, the empirical interpolation method (EIM)~\citep{morBarMNetal04} was introduced to efficiently evaluate non-affinely parametrized functions and quadratic forms in the context of the reduced basis method (RBM)~\citep{morGreetal07}. The EIM relies on a greedy procedure to obtain a suitable `reduced basis' for approximating the function and a stable approach to identify good interpolation points for sparse approximation. The discrete empirical interpolation method (DEIM)~\citep{morChaS10} is a variant of the EIM that makes use of the proper orthogonal decomposition (POD) to obtain a suitable reduced basis to approximate the function. The key distinguishing factor between EIM and DEIM is that while the EIM basis vectors are the (normalized) interpolation error vectors at each iteration of the algorithm, the DEIM basis consists of the POD basis vectors. Furthermore, the EIM selects interpolation points by iteratively approximating the snapshots of the nonlinear function, whereas the DEIM selects the interpolation points by iteratively approximating the POD basis of those snapshots. In this work, we limit our focus to the DEIM.

The DEIM has seen significant success and is widely adopted for approximating functions using sparsely measured data~\citep{Manetal18}. Originally proposed for the approximation of discretized nonlinear functions arising in model order reduction, over the years, it has been found useful in other contexts such as sensor placement~\citep{CleHK21}, interpretable feature selection~\citep{SorE16}, subset selection~\citep{CheFB21}, and data classification~\citep{HenL24}. 

The problem that (D)EIM seeks to address is the following~\citep{morBarMNetal04,morChaS10}. Consider the vector-valued, multivariate nonlinear function denoted by $\bff(\p)$ with $\bff\,:\, \R^{d} \rightarrow \R^{N}$. Suppose that experimental measurements/simulation data of this function, denoted by the set of data snapshots $\Gamma := \{\bff(\p_{1}), \bff(\p_{2}), \ldots, \bff(\p_{m}) \}$, are available at a set of $m$ samples of the argument $\p$, viz., $\Xi := \{\p_{1}, \p_{2}, \ldots, \p_{m}\}$. Using this data, the DEIM algorithm proposes a numerical recipe to construct an approximation of the nonlinear function at any other parameter $\p$, given that data/measurements of $\bff(\p)$ are available at $n$ indices denoted by $\Pi := \{p_{1}, p_{2}, \ldots, p_{n}\}$ with $\Pi \subset \{1, 2, \ldots, N\}$ and $n \ll N$. Concretely, the DEIM algorithm seeks a linear approximation of the nonlinear function, viz.,
\begin{equation}
\label{eq:deim_ansatz}
\begin{aligned}
	\bff \approx \bU \, \bc
\end{aligned}
\end{equation}
where $\bU \in \R^{N \times n}$ is a basis matrix and $\bc \in \R^{n}$ is the vector of unknown coefficients to be determined. Since the problem is over-determined (as $N \gg n$), DEIM imposes an interpolatory projection condition in order to identify the unknown coefficients $\bc$. This reads
\begin{equation}
\label{eq:deim_interp_proj}
\begin{aligned}
\bP^{\tpose} \bff \mathop{=}\limits^! \bP^{\tpose} \bU \, \bc.
\end{aligned}
\end{equation}
The matrix $\bP \in \R^{N \times n}$ is a permutation matrix (consisting only of $1$s and $0$s) that is obtained by selecting a few columns of the identity matrix $\bI \in \R^{N \times N}$. Given that $\bP^{\tpose} \bU$ is invertible, the coefficient $\bc$ is simply
$$
	\bc = \left(\bP^{\tpose} \bU\right)^{-1} \bP^{\tpose} \bff.
$$
Using the expression for the unknown coefficient, the DEIM approximation \cref{eq:deim_ansatz} of the nonlinear function $\bff$ is
\begin{equation}
\label{eq:deim_approx}
\begin{aligned}
\bff \approx \bU \, \left(\bP^{\tpose} \bU\right)^{-1} \bP^{\tpose} \bff.
\end{aligned}
\end{equation}
Since the matrix $\bP$ is a permutation matrix, the term $\bP^{\tpose} \bff$ effectively samples/selects the entries from $n$ rows of the function $\bff$, with the indices of the rows given by the set $\Pi$. This reveals the key benefit of using DEIM for function approximation --- \emph{given the basis matrix $\bU$ and the sampling matrix $\bP$, the nonlinear function at any input $\p$, i.e, $\bff(\p)$ can be reconstructed by simply querying its value at a small set of $n \ll N$ indices}. To ensure the solvability of \cref{eq:deim_interp_proj}, a careful construction of $\bU$ and $\bP$ is required. To construct $\bU$, the singular value decomposition (SVD) of the dataset at the configuration parameters, viz., $\Gamma$ is performed. The matrix $\bU$ consists of the first $n$ columns of the resulting left singular vector matrix. Subsequently, to ensure the invertibility of $\bP^{\tpose} \bU$ in \cref{eq:deim_interp_proj}, a greedy sampling of the rows of $\bU$ is performed to identify the sampling matrix $\bP$; see \citep[Algorithm 1]{morChaS10} for the details. Depending on the application under consideration, the function $\bff$ could represent the states of a large-scale system, or the values of a physical variable measured on some grid. Also, based on the context, the argument $\p$ of the function can be time, system parameters, initial conditions, etc. 

\subsection{Related works}
\label{subsec:intro_related_work}
The problem of reconstructing a function, given its measurements at a set of sparse locations within its domain, is a classical problem in approximation theory and has also been studied in the context of compressed sensing~\citep{Bar07}. In compressed sensing, a universal basis (such as a Fourier basis) is employed and a sparsity-enforcing optimization procedure is used to determine the locations where the function needs to be sampled. In contrast, methods like DEIM use the available measurement data to learn/discover an empirical basis (e.g., using a POD).

DEIM belongs to the family of hyperreduction methods originating from the Gappy-POD~\citep{EveS95} approach. Other methods in this family include the Empirical Interpolation Method (EIM)~\citep{morBarMNetal04}, the Missing Point Estimation (MPE) method~\citep{Astetal08}, the Best Point Interpolation Method (BPIM)~\citep{NguPP08}. All of these methods are closely related and they mainly differ in how the basis matrix is chosen and what problem formulation is adopted to identify the sampling matrix $\bP$. While DEIM, Gappy-POD, and MPE use an orthonormal POD basis to construct $\bU$, the EIM relies on a normalized, non-orthonormal basis obtained directly from the measurement data of the function $\bff(\p)$.

Numerous extensions of DEIM have been proposed in the decade since its introduction, addressing various issues associated with its original formulation and extending its use to a wide range of system classes. We briefly mention some noteworthy works here, while emphasizing that this is not a comprehensive list of all available DEIM variants~\footnote{Given the many variants of the DEIM and its wide usage, the time is opportune for a comprehensive review paper.}. The localized DEIM (LDEIM) approach~\citep{Pehetal14} seeks to improve function approximation for problems evincing a wide range of system behaviours by constructing local DEIM approximants, through an appropriate clustering of the parameter domain. In \citep{PehW15}, the adaptive DEIM (ADEIM) is introduced with the aim of adapting the DEIM basis $\bU$ in a streaming fashion, as new data about the nonlinear function $\bff$ becomes available online. The adaptation of the DEIM basis is done using rank-one updates to the basis matrix. Extensions of DEIM to deal with the preservation of particular problem structure (such as Hamiltonian, port-Hamiltonian, etc.) have also been proposed~\citep{ChaBG16,PagV23}. In~\citep{PehDG20}, an extension of DEIM addressing noisy data is discussed. More recently, a randomized version of DEIM has been proposed~\citep{Sai20} to bring down the potentially large cost of performing SVD of the data matrix. In~\citep{ChoCA20}, the authors focus on reducing the cost of DEIM in the context of model order reduction.
Finally, of particular note is the Q-DEIM approach~\citep{morDrmG16}. It is a variant of DEIM that employs a QR decomposition routine to identify the non-zero (unity) indices of the Boolean sampling matrix $\bP$. The authors show that QR-based approach provably improves the approximation error over standard DEIM. As a result, in practice, the Q-DEIM approach is often the variant of choice instead of the standard DEIM. For this reason, the Q-DEIM will be our variant of choice in this work.

\subsection{Shortcomings of DEIM}
\label{subsec:intro_deim_shortcomings}
The standard Q-DEIM (and also the many DEIM variants mentioned previously) is suitable only when the nonlinear function $\bff(\p)$ is vector-valued. Therefore, the function snapshots are assembled as a matrix with the columns consisting of the vectorized function evaluated at different samples of $\p$. In several applications, the function/data to be approximated is in tensor form (e.g., $\bff(\p) \in \R^{N_{1} \times N_{2} \times N_{3}}$). Instances of this include: spatio-temporal snapshots of parametric partial differential equations (PDEs), sensor data measured at select locations over discrete time steps, for different system configurations. Additionally, in some biological applications, the available experimental dataset is tensor-structured, e.g., brain-machine interface datasets~\citep{Vyasetal18}, or electro-encephalography (EEG) data~\citep{Congetal15}. Matricizing such datasets can distort the data, smearing out any geometric/coherent structures. While there exists the M-DEIM variant that approximates matrix-valued nonlinear functions~\citep{BonMQ17}, it still involves a vectorization of the matrix-valued data in order to apply the DEIM. This poses a serious obstacle when the data is structured, such as 2-D or 3-D data from a mesh.

Very recently, an extension of Q-DEIM, named as higher order DEIM, or HO-DEIM~\citep{Kir22}, was proposed to account for tensor-valued data arising in multi-linear model reduction. The author proposes to construct the DEIM basis through a higher-order singular value decomposition (HO-SVD). More specifically, the matricization/mode-unfolding~\citep{KolB09} of the data tensor is performed and a separate Q-DEIM approximation of each mode unfolding of the data is done. While this addresses the issue of loss of the structure of the data to an extent, nevertheless, the application of an SVD following a matricization of the data can lead to a sub-optimal approximation. Moreover, applying multiple SVDs can be computationally expensive. The present work addresses this issue by proposing an extension of Q-DEIM for tensor-valued data that does not rely on any type of matricization of the tensor data. 

\subsection{Our contributions}
\label{subsec:intro_contributions}
Leveraging the tensor t-product algebra~\citep{KilMP08,KilM11}, the main contribution of this work is an extension of the Q-DEIM to approximate three-dimensional, tensor-valued data. The t-product representation enjoys superior approximation properties over other tensor representations~\citep{Kiletal21} and is therefore the method of choice in this work. Furthermore, it does not involve matricizations of the data tensor, thus preserving relevant geometric information present in the data.

In our proposed approach, hereafter called t-Q-DEIM (short for tensor t-product Q-DEIM), we set up and discuss the solution to an analogous problem as in \cref{eq:deim_ansatz}, but now based on the tensor t-product. We further discuss the equivalent interpolatory projection condition (see \cref{eq:deim_interp_proj}). To identify a suitable basis tensor, we rely on the t-product singular value decomposition (t-SVD). Moreover, to ensure a good choice of the sampling tensor, we develop a strategy (equivalent to the QR decomposition based approach used in Q-DEIM) using the t-product pivoted  QR decomposition (t-pQR). As the central theoretical contribution, we develop a rigorous, computable error bound for the error between the true data and its t-product Q-DEIM-based approximation. An a-priori estimate for one of the terms involved in our error bound is also proposed.
We further illustrate how the derived error bound validates our computational strategy to identify the basis and the sampling indices. 

We note that the proposed t-Q-DEIM is mainly targeting third-order tensor datasets obtained from, e.g., simulations of parametric time-dependent PDEs. As such, the preferred third-order tensor structure for the t-Q-DEIM approximation has the dimension \textsf{space} $\times$ \textsf{time} $\times$ \textsf{parameter} (or \textsf{space} $\times$ \textsf{parameter} $\times$ \textsf{time}). Such third-order tensor datasets arise naturally in PDEs and also in cases where sensor measurements of a physical quantity are available at selected locations, over a set of time samples. If the approximation (based on sparsely measured data) needs to be carried out for many such systems, having different parametric values, then the t-Q-DEIM would be ideally suited for this use case. It is worth noting that such third-order data tensors have been considered also in the case of parametric DMD approximations~\citep{pdmd_huhn23,Andreuzzi23,morSunFCetal23}. Notwithstanding this initial motivation, we have observed that the proposed t-Q-DEIM approach also performs well for other third-order tensor data, such as those common in biology. In this context, one of the numerical examples considered in this work deals with experimental data of the form \textsf{neurons} $\times$ \textsf{time} $\times$ \textsf{trials}. 

\subsection{Organization}
\label{subsec:intro_organization}
In \Cref{sec:mathprelim}, we discuss the mathematical background of the tensor t-product, highlighting only those aspects relating to the t-Q-DEIM. The main contributions of this work are presented in \Cref{sec:tqdeim}. We introduce the problem formulation of t-Q-DEIM and discuss the interpolatory projection condition. We further discuss the computation of the key quantities involved. The main theorem and the associated error analysis are also presented. \Cref{sec:numerics} relates to the numerical experiments. We validate t-Q-DEIM on five numerical examples, representing a wide variety of systems. We compare its performance with the standard Q-DEIM, emphasizing the computational benefits that t-Q-DEIM offers while also providing a better approximation of the target data. We also numerically compare our proposed t-pQR based sampling method with a closely related existing approach. We conclude in~\Cref{sec:conclusion} offering a summary of the current work and highlighting some potential areas for future research.

\section{Mathematical preliminaries}
\label{sec:mathprelim}
We briefly review below some mathematical background in support of the main contributions to be presented in~\Cref{sec:tqdeim}. Before that, we highlight the mathematical notations to be used throughout the discussion.

\subsection{Notations}
Throughout this work
\begin{itemize}
	\item bold alphabets in lower-case letters shall denote vectors, e.g., $\mathbf{a} \in \R^{N}$,
	\item bold alphabets in upper-case letters shall denote matrices, e.g., $\mathbf{A} \in \R^{N \times n}$,
	\item upper-case Greek alphabets shall denote sets, e.g., $\Gamma$,
	\item calligraphic alphabets in upper-case letters shall denote third-order tensors, e.g., $\cA \in \R^{N \times p \times N}$,
	
	\item calligraphic alphabets in upper-case letters and a cap on top denote the Fourier domain representation of the third-order tensor, e.g., $\chA$,
	
	\item bold lower-case Fraktur letter, e.g., $\mathfrak{f}$ denotes a slice, along the second dimension of a third order tensor, e.g., $\mathfrak{f} = \cF(:,j,:) \in \R^{N \times 1 \times M}$, $\cF \in \R^{N \times n \times M}$,
	
	\item script-style upper-case letters, e.g. $\mathscr{D}$, are used to denote an operator,
	
	\item we use \MATLAB notation to denote row/columns of a matrices/tensors, e.g., $\mathbf{A}(:, i)$ refers to the $i$-th column of the matrix $\mathbf{A}$ and $\mathbf{A}(j, :)$ refers to the $j$-th row of the matrix $\mathbf{A}$,
	\item unless specified otherwise, the lower-case alphabets, $i, j, k$ are reserved to be used as indices,
	
	\item we denote by blackboard upper-case letters $\R^{}$ and $\C^{}$ the set of real, complex numbers,  respectively, and use $\mathbb{K}$ to mean a general field; by extension, the same symbols with appropriate superscript dimensions shall denote the space of matrices/tensors, e.g., $\R^{N \times n}$ is the set of real-valued matrices having dimension $N \times n$.
\end{itemize}

\subsection{The t-product algebra}
The tensor t-product was first introduced in~\citep{KilMP08}, specifically addressing the extension of the matrix SVD to third-order tensors. The theoretical underpinnings were further developed in~\citep{Bra10, KilM11}. In what follows, we only discuss the elements of the t-product algebra relevant for the problem under consideration and refer the interested reader to the work~\citep{Kiletal21} for a more complete discussion.

Consider a third-order tensor $\cA \in \R^{m \times \ell \times q}$.  Being a third-order tensor, three different orientation of the tensor are conceivable. We detail these below.
\begin{definition}[Horizontal slices]
	The horizontal slices of the third-order tensor $\cA$ refer to the $m$ matrices formed by fixing the first index of $\cA$, viz., $\cA\left(i,\,:,\,:\right) \in \R^{\ell \times q}$, $i =1, 2, \ldots, m$.
\end{definition}

\begin{definition}[Lateral slices]
	The lateral slices of the third-order tensor $\cA$ refer to the $\ell$ matrices formed by fixing the second index of $\cA$, viz., $\cA\left(:,\,j,\,:\right) \in \R^{m \times q}$, $j =1, 2, \ldots, \ell$.
\end{definition}

\begin{definition}[Frontal slices]
	The frontal slices of the third-order tensor $\cA$ refer to the $q$ matrices formed by fixing the third index of $\cA$, viz., $\cA\left(:,\,:,\,k\right) \in \R^{m \times \ell}$, $k =1, 2, \ldots, q$.
\end{definition}
Above and also in the discussion to follow, we have denoted third-order tensors as being real-valued. We emphasize that they can also be complex-valued.
We next define the tube fiber of a third-order tensor. 
\begin{definition}[Tube fiber~\citep{Kiletal13}]
	Let $\cA \in \R^{m \times \ell \times q}$ be a third-order tensor. We define $\mathbf{a}_{ij} := \cA(i,\, j,\, :) \in \R^{1 \times 1 \times q}$  as the $ij$-th tube fiber of the tensor $\cA$ having tubal length $q$.	
\end{definition}

A tube fiber is the t-product algebra equivalent of a scalar in $\mathbb{R}$. For this reason, we also refer to a tube fiber as a tubal scalar~\footnote{We will use the term \emph{tube fiber} when we refer to a tubal scalar of a particular tensor.}. The space of all tubal scalars having $q$ entries is denoted $\K_{q}$. Based on the previous definition, a third-order tensor $\cB \in \R^{m \times 1 \times q}$ is a length $m$ vector of tubal scalars and is said to belong to the space $\K_{q}^{m}$. Further, a third-order tensor $\cB \in \R^{m \times \ell \times q}$ is an $m \times \ell$ matrix of tubal scalars and belongs to the space $\K_{q}^{m \times \ell}$; see~\citep{Bra10} for a detailed discussion. 

\begin{definition}[Block circulant matrix]
	A block circulant matrix of the third-order tensor $\cA \in \R^{m \times \ell \times q}$ is the circulant matrix formed using the frontal slices of $\cA$, viz., $\bA^{(k)} := \cA\left(:,\,:,\,k\right)$, $k = 1, 2, \ldots, q$. It is given by
	$$
		\bcirc{A} := 
		\begin{bmatrix}
			\bA^{(1)} & \bA^{(q)} & \cdots & \bA^{(2)}\\
			\bA^{(2)} & \bA^{(1)} & \cdots & \bA^{(3)}\\
			\vdots    & \vdots    &  \ddots & \vdots\\[0.3em]
			\bA^{(q)} & \bA^{(q-1)} & \cdots & \bA^{(1)}
		\end{bmatrix} \in \R^{mq \times \ell q}.
	$$
\end{definition}

\begin{definition}[Unfolding operation]
	The unfolding of the third-order tensor $\cA$ is defined as the operation that vertically stacks the frontal slices of $\cA$. That is,
	$$
		\unfold{A} := 
		\begin{bmatrix}
			\bA^{(1)}\\
			\bA^{(2)}\\
			\vdots\\[0.3em]
			\bA^{(q)}
		\end{bmatrix} \in \R^{mq \times \ell}.
	$$
\end{definition}

\begin{definition}[Folding operation]
	The folding operation is the inverse of the unfolding operation. The fold operation stacks, along the third dimension, the blocks of an unfolded third-order tensor. That is,
	$$
		\fold{\unfold{A}} := \cA.
	$$
\end{definition}
Relying on the above definitions, the third-order tensor product (t-product) of two tensors is defined as follows.

\begin{definition}[t-product~\citep{KilMP08}]
	Consider the third-order tensors $\cA \in \R^{m \times \ell \times q}$ and $\cB \in \R^{\ell \times r \times q}$. The third-order tensor product (t-product) of $\cA$, $\cB$ denoted \tprod{A}{B} is the $m \times r \times q$ tensor
	$$
		\tprod{A}{B} = \fold{\bcirc{A} \cdot \unfold{B}}.
	$$
\end{definition}
Note that, for the above definition to be valid, the second dimension of the tensor $\cA$ and the first dimension of the tensor $\cB$ need to be the same. Moreover, the third dimensions of both tensors, viz. $q$, should match.

\begin{remark}
	\label{rem:fft_multiply}
	The t-product, as it is defined, has close connections to the Fourier transform. This connection can be leveraged to efficiently evaluate the t-product for dense third-order tensors. The key steps are as follows~\citep{KilM11}:
	\begin{enumerate}
		\item To evaluate the t-product $\tprod{A}{B}$ where $\cA \in \R^{m \times \ell \times q}$ and $\cA \in \R^{\ell \times r \times q}$, compute their respective Fast Fourier Transform (FFT) denoted by $\chA :=\textnormal{ \texttt{fft}}(\cA) \in \C^{m \times \ell \times q}$ and $\chB :=\textnormal{ \texttt{fft}}(\cB) \in \C^{\ell \times r \times q}$, with the FFT being applied along the third dimension. In \MATLAB notation, this reads $\chA :=\textnormal{ \texttt{fft}}(\cA,\, [\,],\, 3)$.
		
		\item Perform the pairwise matrix multiplication of each of the frontal slices of $\chA, \chB$, i.e., 
		$$
			\chC\left(:,\,:,\,k\right) := \chA\left(:,\,:,\,k\right) \cdot  \chB\left(:,\,:,\,k\right)
		$$
		for $k = 1, 2, \ldots, q$. For future use, we refer to this frontal slice-wise matrix multiplication of two third-order tensors with the notation 
		$$
			\chC := \chA\,\, \triangle\,\, \chB.
		$$
		\item The t-product is then simply the inverse FFT of the tensor $\chC$, i.e.,
		$$
			\cC := \tprod{A}{B} = \textnormal{\texttt{ifft}}(\chC).
		$$
	\end{enumerate} 

The cost of computing the t-product via the FFT is at most $\mathcal{O}(m \ell r q)$ flops~\textnormal{\citep{KilM11}}.
\end{remark}

An important benefit of using the t-product framework is that, it offers a generalization to the tensor setting of key linear algebra concepts such as the identity matrix, the inverse of a matrix, the transpose, and orthogonality. We briefly define these concepts to support the definition of the tensor SVD (t-SVD).

We define t-linearity generalizing the notion of linearity in the vector space to the t-product algebra.
\begin{definition}[t-linearity~\citep{Kiletal13}]
	Let $\{\mathbf{b}_{i}\}_{i=1}^{k}$ be $k$ tubal scalars with $\mathbf{b}_{i} \in \K_{q}$. A t-linear combination of the tensors $\cA_{j} \in \K_{q}^{m}$, $j = 1, 2, \ldots, k$ is defined as
	$$
	\cA_{1} * \mathbf{b}_{1} + \cA_{2} * \mathbf{b}_{2} + \cA_{3} * \mathbf{b}_{3} + \cdots + \cA_{k} * \mathbf{b}_{k}.
	$$
\end{definition}

\begin{definition}[t-transpose~\citep{KilM11,Kiletal13}]
	The transpose of the third-order tensor $\cA \in \R^{m \times \ell \times q}$ is defined as the $\ell \times m \times q$ tensor obtained by performing the matrix transpose of each of its frontal slices, viz., $(\bA^{(k)})^{T}$, $k = 1, 2, \ldots, q$ followed by reversing the order of the transposed frontal slices from $2$ through $q$. We have
	$$
		\cA^{\textnormal{\tpose}} := \textnormal{\texttt{fold}}\left(
		\begin{bmatrix}
			\left(\bA^{(1)}\right)^{\textnormal{\tpose}}\\
			\left(\bA^{(q)}\right)^{\textnormal{\tpose}}\\
			\left(\bA^{(q-1)}\right)^{\textnormal{\tpose}}\\
			\vdots\\[0.3em]
			\left(\bA^{(2)}\right)^{\textnormal{\tpose}}
		\end{bmatrix} 
		\right) \in \R^{\ell \times m \times q}.
	$$
\end{definition}

\begin{definition}[t-identity~\citep{KilM11}] 
	We define the third-order identity tensor \cIdim{m}{q} $\in \R^{m \times m \times q}$ as the tensor whose first frontal slice is the $m \times m$ identity matrix and whose remaining frontal slices are zero matrices.
\end{definition}

\begin{definition}[Permutation tensor~\citep{KilMP08,Haoetal13}]
	A third-order tensor $\cP \in \R^{m \times m \times q}$ is a permutation tensor if its entries consist only of zeros and ones and furthermore it has exactly $m$ entries of unity and such that if $\cP(i,\,j,\,k) = 1$, it is the only non-zero entry in the $i$-th row, $j$-th column, and $k$-th slice, where the term slice denotes the third dimension. Additionally, 
	$$
		\cP^{\textnormal{\tpose}} * \cP = \cP * \cP^{\textnormal{\tpose}} = \cIdim{m}{q}.
	$$
\end{definition}

\begin{definition}[t-orthogonality~\citep{KilM11}]
	\label{def:t-orth}
	The third-order tensor $\cA \in \R^{m \times m \times q}$ is orthogonal if 
	$$
		\cA^{\textnormal{\tpose}} *  \cA = \cA * \cA^{\textnormal{\tpose}} = \cI.
	$$
\end{definition}
The orthogonality can be interpreted in terms of each lateral slice of  the orthogonal matrix $\cA$. It holds that
\begin{align}
	\label{eq:t-orth-property}
	\cA(:,\,i,\,:)^{\tpose} * \cA(:,\,j,\,:) = \begin{cases}
	\mathfrak{i}_{1 \times 1}, & \text{if $i = j$}.\\
	\mathfrak{o}_{1 \times 1}, & \text{otherwise}.
	\end{cases}
\end{align}
In the above, $\mathfrak{i}_{1 \times 1} \in \R^{1 \times 1 \times M}$ denotes a tubal fiber with $1$ in the first frontal slice and zeros everywhere else while $\mathfrak{o}_{1 \times 1} \in \R^{1 \times 1 \times M}$ is a zero tubal fiber with zeros everywhere.

\begin{definition}[t-inverse~\citep{KilM11}]
	The third-order tensor $\cA \in \R^{m \times m \times q}$ is said to be invertible if there exists a tensor $\cB \in \R^{m \times m \times q}$ such that
	$$
		\tprod{A}{B} = \cIdim{m}{q}
	$$ and
	$$
		\tprod{B}{A} = \cIdim{m}{q}.
	$$
\end{definition}
The notion of norm in the t-product algebra, particularly the Frobenius norm, is analogous to the matrix case.
\begin{definition}[t-Frobenius norm]
	The Frobenius norm of the third-order tensor $\cA \in \R^{m \times \ell \times q}$ is defined as
	$$
	\norm{\cA}{F} = \sqrt{\sum_{i=1}^{m} \sum_{j=1}^{\ell} \sum_{k=1}^{q} a_{ijk}^{2}}
	$$
	where $a_{ijk} \in \R^{}$ is the element in $\cA(i,\,j,\,k)$.
\end{definition}
An extension of the spectral norm to the t-product algebra was introduced in~\citep{Luetal20} and is defined next.
\begin{definition}[t-spectral norm]
	\label{def:tsn}
	The t-spectral norm $\norm{\cdot}{}$ of a third-order tensor $\cA \in \R^{m \times \ell \times q}$
	is defined as
	$$
		\norm{\cA}{} := \norm{\bcirc{A}}{2} = \norm{\chA}{}
	$$
	where $\bcirc{A} \in \R^{mq \times \ell q}$ and $\chA :=\textnormal{ \texttt{fft}}(\cA)$.
\end{definition}
The t-spectral norm can be computed by transforming $\cA$ into its Fourier domain equivalent $\chA$ and taking the maximum of the matrix spectral norm over all frontal slices of $\chA$. In essence, 
$$
	\norm{\cA}{} = \norm{\chA}{} = \max \limits_{i = 1, 2, \ldots, q} \norm{\chA(:,\,:,\,i)}{2}.
$$
\begin{theorem}[t-SVD~\citep{KilM11}]
	Consider the third-order tensor $\cA \in \R^{m \times \ell \times q}$. $\cA$ can be factorized as
	$$
		\cA = \cU * \cS * \cW^{\textnormal{\tpose}}
	$$
	with $\cU \in \R^{m \times m \times q}$ and $\cW \in \R^{\ell \times \ell \times q}$ being orthogonal tensors and $\cS \in \R^{m \times \ell \times q}$ is an \emph{f-diagonal} tensor, meaning each of its frontal slices is a diagonal matrix.
\end{theorem}
The t-SVD enjoys an optimal approximation property, analogous to the matrix SVD. The following theorem illustrates this.
\begin{theorem}[t-SVD optimality~\citep{KilM11}]
	Consider the third-order tensor $\cA \in \R^{m \times \ell \times q}$ and its t-SVD given by $\cA = \cU * \cS * \cW^{\textnormal{\tpose}}$. Let $\cA_{n}$ represent the $n$-term t-SVD approximation with
	$$	
		\cA_{n} = \sum_{i=1}^{n} \cU(:,\,i,\,:) * \cS(i,\,i,\,:) * \cW(:,\,i,\,:)^{\textnormal{\tpose}}
	$$ 
	where $n < \min(m, \ell)$.	Then, it holds that 
	$$
		\cA_{n} = \arg \min \limits_{\tilde{\cA}}	\norm{\cA - \tilde{\cA}}{F}	
	$$
\end{theorem}
The above optimal approximation property theorem for the t-SVD is analogous to the famous \textit{Schmidt-Eckhart-Young-Mirsky theorem} for the matrix SVD. Next, we define the t-product versions of the matrix QR and pivoted QR decomposition.
\begin{definition}[t-QR decomposition~\citep{KilMP08}]
	Consider the third-order tensor $\cA \in \R^{m \times \ell \times q}$. The tensor $\cA$ can be factored as
	$$
		\cA = \cQ * \cR
	$$
	where $\cQ \in \R^{m \times m \times q}$ is an orthogonal tensor and $\cR \in \R^{m \times \ell \times q}$ is an f-upper triangular tensor, where f-upper triangular means that each of the frontal slices of $\cR$ is an upper triangular matrix.
\end{definition}
\begin{definition}[Pivoted t-QR decomposition~\citep{Haoetal13}]
	Consider the third-order tensor $\cA \in \R^{m \times \ell \times q}$. The tensor $\cA$ can be factored as
	$$
		\cA * \cP = \cQ * \cR
	$$
	where $\cQ \in \R^{m \times m \times q}$ is an orthogonal tensor, $\cR \in \R^{m \times \ell \times q}$ is an f-upper triangular tensor and $\cP \in \R^{\ell \times \ell \times q}$ is a permutation tensor.
\end{definition}
The computation of the pivoted t-QR decomposition is detailed in~\Cref{appendix:tpqr}.


\section{The tensor t-product Q-DEIM (t-Q-DEIM)}
\label{sec:tqdeim}
The section presents and develops the main contributions of this work, viz., a tensor t-product extension of the standard DEIM and the corresponding error analysis. We start with the description of the problem setup.
\subsection{t-Q-DEIM problem setup}
We are interested in the approximation of a nonlinear, tensor-valued function\\ $\mathfrak{f}(\p)\,:\, \R^{d} \rightarrow \R^{N \times 1 \times M}$, based on sparsely measured data. We seek a t-linear approximation of this function as below
\begin{equation}
\label{eq:tpdeim_ansatz}
	\begin{aligned}
		\mathfrak{f} \approx \left(\cU_{1} * \mathbf{c}_{1} + \cU_{2} * \mathbf{c}_{2} + \cdots + \cU_{n} * \mathbf{c}_{n}\right) = \cU * \mathfrak{c}
	\end{aligned}
\end{equation}
where $\cU := \left[\cU_{1}, \cU_{2}, \dots, \cU_{n}\right] \in \R^{N \times n \times M}$ is the basis tensor with $n \ll N$. Note that $\cU_{i} \in \R^{N \times 1 \times M}$ is a lateral slice of $\cU$. Further, $\mathfrak{c} \in \R^{n \times 1 \times M}$ is the tensor of unknown coefficients and let
$$
	\mathfrak{c} :=
	\begin{bmatrix}
		\mathbf{c}_{1}\\
		\vdots\\
		\mathbf{c}_{n}
	\end{bmatrix}
$$
with $\mathbf{c}_{i} \in \R^{1 \times 1 \times M}$ being a horizontal slice of $\mathfrak{c}$. The above ansatz to approximate $\mathfrak{f}$ will be accurate as long as the dominant information of the tensor-valued function lies in the \emph{t-linear} span of the lateral slices of the basis matrix $\cU$~\citep{Kiletal13}. We will ensure this through an appropriate construction of the basis tensor. Evidently, the above problem is over-determined, as there are $N M$ equations and $n M$ unknowns. We enforce an interpolatory projection to make the solution amenable. To this end, we require
\begin{equation}
\label{eq:tpdeim_interp_proj}
\begin{aligned}
	\cP^{\tpose}\, * \mathfrak{f}\, \mathop{=}\limits^!\, \cP^{\tpose} *\, \cU\, *\, \mathfrak{c}
\end{aligned}
\end{equation}
where $\cP \in \R^{N \times n \times M}$ is a sampling tensor that picks/samples exactly $n$ horizontal slices of $\mathfrak{f}$. Notice that, the system is no longer over-determined as now there are exactly $nM$ equations and unknowns. Assuming for a moment that  $\left(\cP^{\tpose} *\, \cU\right)^{-1}$ is invertible, the unknown coefficients can be expressed as
\begin{equation}
\label{eq:tpdeim_coeff}
\begin{aligned}
	\mathfrak{c} = \left( \cP^{\tpose} *\, \cU\right)^{-1}\, *\, \cP^{\tpose}\, *\, \mathfrak{f}.
\end{aligned}
\end{equation}
Substituting \cref{eq:tpdeim_coeff} into \cref{eq:tpdeim_ansatz}, we obtain the  t-Q-DEIM approximation $\mathfrak{f}_{\text{tq}}(\p)$ of the nonlinear function $\mathfrak{f}(\p)$:
\begin{equation}
\label{eq:tpdeim_approx}
\begin{aligned}
	\mathfrak{f} \approx \mathfrak{f}_{\text{tq}} := \cU * \left( \cP^{\tpose} * \cU\right)^{-1} * \cP^{\tpose} * \mathfrak{f}.
\end{aligned}
\end{equation}
Observe that once $\cU$ and $\cP$ are known, the quantity $\left(\cU * \left( \cP^{\tpose} * \cU\right)^{-1}\right) \in \R^{N \times n \times M}$ can be precomputed and stored. An evaluation of the t-Q-DEIM approximation at a new parameter will simply consist of a t-product between this precomputed quantity and the function evaluated at $n$ horizontal slices. The cost of this is at most $\mathcal{O}(nNM)$. 

Let $\mathscr{D} := \cU*(\cP^{\tpose}* \cU)^{-1}*\cP^{\tpose} \in \R^{N \times N \times M}$ denote the third-order t-Q-DEIM projection operator which acts on $\mathfrak{f}$ to produce the t-Q-DEIM approximation. Using this, the t-Q-DEIM approximation can be interpreted as an \emph{interpolatory projection} of the function $\mathfrak{f}$ on to the t-span of the lateral slices of the projection operator. For the projection operator $\mathscr{D}$, the following properties hold:
\begin{property}[Projection]
	\label{prop:t-proj}
	$$
	\mathscr{D}^{2} = \mathscr{D}*\mathscr{D} = \left(\cU*(\cP^{\tpose}* \cU)^{-1}*(\cP^{\tpose}*\cU)*(\cP^{\tpose}* \cU)^{-1}*\cP^{\tpose}\right) = \left(\cU*(\cP^{\tpose}* \cU)^{-1}*\cP^{\tpose}\right) = \mathscr{D}.
	$$
\end{property}
Note that we have used the associative property of the t-product~\citep{Kiletal13} to obtain the above result.
\begin{property}[Interpolatory property]
	\label{prop:t-interp_proj}
	$$
	\mathcal{P}^{\tpose} * (\mathscr{D}*\mathfrak{f}) = (\mathcal{P}^{\tpose} * \cU)*(\cP^{\tpose}* \cU)^{-1}*\cP^{\tpose}*\mathfrak{f} = \cP^{\tpose}*\mathfrak{f}.
	$$
\end{property}
In essence, \Cref{prop:t-interp_proj} states that the t-Q-DEIM approximant \emph{interpolates} the original function $\mathfrak{f}$ at the indices encoded in the sampling tensor. Both the above properties rely on the assumption that $(\mathcal{P}^{\tpose} * \cU)$ is invertible. Next, we discuss the computation of the basis tensor $\cU$ and the corresponding sampling tensor $\cP$.

\subsection{Computing the t-Q-DEIM basis and sampling tensor}
To compute the tensor basis $\cU$ and sampling tensor $\cP$ in the t-Q-DEIM approximation \cref{eq:tpdeim_approx}, we leverage the t-product algebra tools discussed in \Cref{sec:mathprelim}.

\paragraph{t-SVD to compute the basis}
We suppose that experimental measurement data/artificial simulation data of the nonlinear function $\mathfrak{f}(\p)$ is available at a set of $n_{s}$ samples of $\p$, i.e.,\\ $\Xi := \{\p_{1}, \p_{2}, \ldots, \p_{n_{s}}\}$. Let $\cF \in \R^{N \times n_{s} \times M}$ be the third-order tensor whose $j$-th lateral slice contains the evaluation of the nonlinear function $\mathfrak{f}(\p_{j})$ at the $j$-th parameter, with $j = 1, 2, \ldots, n_{s}$. We proceed by using the t-SVD to write the data tensor as
\begin{align}
	\label{eq:snapshot_svd}
	\cF = \cU_{f} * \cS_{f} * \cW_{f}^{\tpose}
\end{align}
with $\cU_{f} \in \R^{N \times N \times M}$, $\cW_{f} \in \R^{n_{s} \times n_{s} \times M}$ being orthogonal tensors and where $\cS_{f} \in \R^{N \times n_{s} \times M}$ is an f-diagonal tensor, each of whose frontal slices is a diagonal matrix. We define the basis tensor $\cU$ as the tensor obtained by choosing the first $n$ lateral slices of the left singular tensor $\cU_{f}$, i.e., $\cU := \cU_{f}\left(:,\, 1:n,\, :\right) \in \R^{N \times n \times M}$. Choosing the basis tensor $\cU$ based on the t-SVD not only enjoys the optimal approximation property, but also offers a well-conditioned basis owing to $\cU_{f}$ being orthogonal.

\paragraph{Pivoted t-QR to compute the sampling tensor}
Having identified a suitable basis tensor by leveraging the t-SVD, we now discuss the computation of an appropriate sampling tensor $\cP$. The role of the sampling tensor is to select a small number of horizontal slices of the nonlinear function $\mathfrak{f}(\p)$. Furthermore, the choice of $\cP$ should also ensure that the term $\left(\cP^{\tpose} * \cU\right)$ in \cref{eq:tpdeim_approx} is invertible. Recall that in \citep{morDrmG16}, a pivoted QR decomposition of the basis matrix $\bU$ was performed to identify the sampling matrix $\bP$. The motivation there was that choosing the rows of $\bU$ corresponding to the pivots indices (of $\bU^{\tpose}$) leads to a better conditioning of the quantity $\bP^{\tpose} \bU$. We proceed in a similar spirit by considering the pivoted t-QR decomposition of the basis $\cU$. First, we obtain the frequency domain representation of $\cU$ by performing its FFT, yielding $\widehat{\cU} := \textnormal{\texttt{fft}}(\cU)$. Next, we extract the first frontal slice of $\widehat{\cU}$, i.e., $\widehat{\mathbf{U}}^{(1)} := \widehat{\cU}(:,\,:,\,1) \in \C^{N \times n}$. Finally, a standard pivoted QR decomposition of $\left(\widehat{\mathbf{U}}^{(1)}\right)^{*}$ is done to obtain the pivots. Then, retaining only the first $n$ pivots we construct the index set $\Pi_{f} := \{p_{1}, p_{2}, \ldots, p_{n}\}$ where $\Pi_{f} \subset \{1, 2, 3, \ldots, N\}$. We define the tensor $\widehat{\mathcal{P}}$ by populating all $M$ of its frontal slices with the same matrix $\mathbf{I}_{s} \in \R^{N \times n}$. The matrix $\mathbf{I}_{s}$ has as its $j$-th column, the $p_{j}$-th column of the standard identity matrix $\mathbf{I} \in \R^{N \times N}$, and $p_{j} \in \Pi_{f}$. In essence, the matrix $\mathbf{I}_{s}$ is simply the permutation of the columns of the standard identity matrix, with the permutation indices given by the set $\Pi_{f}$. Finally, the sampling tensor $\cP$ is obtained after an inverse FFT as $\cP := \textnormal{\texttt{ifft}}(\widehat{\mathcal{P}})$. Ultimately, $\cP$ is a permutation tensor, therefore, its first frontal slice will consist of $\mathbf{I}_{s}$ while the remaining frontal slices consist of zero matrices of appropriate dimension. The pseudocode to compute the sampling indices is sketched in~\Cref{alg:tpqr}.
\begin{remark}
	Once the sampling indices $\Pi_{f}$ are available, computationally speaking, it is not necessary to construct the sampling tensor $\cP$. The horizontal slices can be sampled directly by choosing those rows of the function $\mathfrak{f}$ (or rows of the basis $\cU$) as $\mathfrak{f}(\Pi_{f},\,:,\,:)$ as it holds that
	$$
		\cP^{\tpose} * \mathfrak{f} = \mathfrak{f}(\Pi_{f},\,:,\,:)
	$$ and, similarly,
	$$
		\cP^{\tpose} * \cU = \cU(\Pi_{f},\,:,\,:).
	$$
\end{remark}
\begin{algorithm}[t!]
	\SetKwInOut{Input}{Input}
	\SetKwInOut{Output}{Output}
	\Input{Basis tensor $\cU \in \R^{N \times n \times M}$\vspace{0.5em}}  
	\Output{Sampling indices $\mathbf{p} \in \R^{n}$\vspace{0.5em}}
	
	{Perform a FFT of $\cU$ along the third dimension to obtain $\chU$\vspace{0.5em}}
	
	{Extract the first frontal slice of $\chU$; $\widehat{\mathbf{U}}^{(1)} := \chU(:,\,:,\,1)$\vspace{0.5em}}
	
	{Perform the pivoted QR decomposition of $\left(\widehat{\mathbf{U}}^{(1)}\right)^{*} \in \R^{n \times N}$ and select the first $n$ pivots, setting $\mathbf{p} := [\mathrm{p}_{1}, \mathrm{p}_{2}, \ldots, \mathrm{p}_{n}]^{\tpose} \in \R^{n}$\vspace{0.5em}}
	\caption{t-pQR}
	\label{alg:tpqr}
\end{algorithm}

We summarize the t-Q-DEIM approach in~\Cref{alg:tqdeim}, where one of the outputs, the tensor $\cD$ is computed only once and stored. It can then be repeatedly used for the evaluation of the interpolation $\mathfrak{f}_{\text{tq}}$~\cref{eq:tpdeim_approx} at any value of $\p$. Furthermore, for comparison, the Q-DEIM is summarized in~\Cref{alg:qdeim}.
\begin{algorithm}[t!]
	\SetKwInOut{Input}{Input}
	\SetKwInOut{Output}{Output}
	\Input{Training data $\cFtrain \in \R^{m \times N_{\text{train}} \times q}$, reduced dimension $n$\vspace{0.5em}}  
	\Output{t-Q-DEIM tensor $\mathcal{D} \in \R^{m \times n \times q}$, sampling locations $\{\mathrm{p}_{1}, \mathrm{p}_{2}, \ldots, \mathrm{p}_{n}\}$\vspace{0.5em}}
	
	{Perform t-SVD on $\cFtrain$ to obtain basis matrix $\cU \in \R^{m \times n \times q}$\vspace{0.5em}}
	
	{Obtain pivot indices $\mathbf{p} := [\mathrm{p}_{1}, \mathrm{p}_{2}, \ldots, \mathrm{p}_{n}]^{\tpose} \in \R^{n}$ through a pivoted t-QR decomposition of the first frontal slice in the Fourier domain representation of $\cU$; see~\Cref{alg:tpqr}\vspace{0.5em}}
	
	{Compute $\mathcal{D} := \cU * \left(\cU(\mathbf{p},\,:,\,:) \right)^{-1} \in \R^{m \times n \times q}$\vspace{0.5em}}
	\caption{t-Q-DEIM}
	\label{alg:tqdeim}
\end{algorithm}
\begin{algorithm}[t!]
	\SetKwInOut{Input}{Input}
	\SetKwInOut{Output}{Output}
	\Input{Training data $\cFtrain \in \R^{m \times N_{\text{train}} \times q}$, reduced dimension $n$\vspace{0.5em}}  
	\Output{Q-DEIM matrix $\mathbf{D} \in \R^{m \times n}$, sampling locations $\{\mathrm{p}_{1}, \mathrm{p}_{2}, \ldots, \mathrm{p}_{n}\}$\vspace{0.5em}}
	{Vectorize input snapshot tensor $\mathbf{F} = \text{vec}(\cFtrain) \in \R^{m \times qN_{\text{train}}}$\vspace{0.5em}}
	
	{Perform SVD on $\mathbf{F}$ to obtain basis matrix $\bU \in \R^{m \times n}$\vspace{0.5em}}
	
	{Obtain pivot indices $\mathbf{p} := [\mathrm{p}_{1}, \mathrm{p}_{2}, \ldots, \mathrm{p}_{n}]^{\tpose} \in \R^{n}$ through pivoted QR decomposition of $\bU^{\tpose}$\vspace{0.5em}}
	
	{Compute $\mathbf{D} := \bU\cdot \left(\bU(\mathbf{p},\,:) \right)^{-1} \in \R^{m \times n}$\vspace{0.5em}}
	\caption{Q-DEIM~\citep{morDrmG16}}
	\label{alg:qdeim}	
\end{algorithm}
Having discussed the t-Q-DEIM approximation of a nonlinear tensor-valued function and the efficient computation of the quantities involved in the approximation, we perform the error analysis of the t-Q-DEIM scheme.

\subsection{Error analysis}
This section derives an error bound for the t-Q-DEIM approximation. Before discussing the main theorem, we first state and prove some lemmas which will be used in the proof of the main theorem.

\begin{lemma}
	\label{lem:t-submultiplicative}
	The t-spectral norm is sub-multiplicative.
	\begin{proof}
		Let $\cA \in \R^{m \times \ell \times q}$ and $\cB \in \R^{\ell \times r \times q}$ be third-order tensors. Using \Cref{def:tsn}, we have
		$$
			\norm{\cA*\cB}{} = \norm{\bcirc{\cA*\cB}}{}.
		$$
		The block circulant matrix of $(\cA * \cB)$ can be block diagonalized~\citep{KilM11} using the normalized discrete Fourier transform matrix $\mathbf{F}_{q} \in \R^{q \times q}$. This yields
		\begin{align*}
			\norm{\cA*\cB}{} &= \norm{\bcirc{\cA*\cB}}{} = \norm{(\mathbf{F}_{q} \otimes \mathbf{I}_{m \times m}) \cdot \bcirc{\cA*\cB} \cdot (\mathbf{F}_{q}^{*} \otimes \mathbf{I}_{r \times r})}{},\\[1em]
			&= \left\lVert 
			\begin{pmatrix}
				\chC_{1} \\
				& \ddots \\
				& & \chC_{q}
			\end{pmatrix}  
			\right\rVert,\\[1em]
			&= \left\lVert 
			\begin{pmatrix}
				\chA_{1} \\
				& \ddots \\
				& & \chA_{q}
			\end{pmatrix}  
			\begin{pmatrix}
				\chB_{1} \\
				& \ddots \\
				& & \chB_{q}
			\end{pmatrix}			
			\right\rVert
		\end{align*}
		where $\chC_{i} = \chA_{i} \chB_{i}$, $i = 1, 2, \ldots, q$ and $\chA, \chB$ are the representations, respectively, of $\cA, \cB$ in the Fourier domain obtained by taking the Fourier transform along the third dimension of the tensors. The second equality above results from the observation that the spectral norm of the block diagonal matrix consisting of the frontal slices in the Fourier domain is the same as that of the block circulant matrix, as the former is obtained through a \emph{unitary transformation} of the latter.	Using the property of the matrix spectral norm, we get
		\begin{align*}
			\norm{\cA*\cB}{} &= \left\lVert 
			\begin{pmatrix}
			\chA_{1} \\
			& \ddots \\
			& & \chA_{q}
			\end{pmatrix}  
			\begin{pmatrix}
			\chB_{1} \\
			& \ddots \\
			& & \chB_{q}
			\end{pmatrix}			
			\right\rVert,\\[1em]
			&\leq \left\lVert 
			\begin{pmatrix}
				\chA_{1} \\
				& \ddots \\
				& & \chA_{q}
			\end{pmatrix}
			\right\rVert
			\left\lVert 
			\begin{pmatrix}
			\chB_{1} \\
			& \ddots \\
			& & \chB_{q}
			\end{pmatrix}
			\right\rVert,\\
			&= \norm{\bcirc{\cA}}{} \norm{\bcirc{\cB}}{},\\
			&= \norm{\cA}{} \norm{\cB}{}.
		\end{align*}
	\end{proof}
\end{lemma}
\begin{lemma}
	\label{lem:tprod_triple}
	For three third-order tensors $\cA \in \R^{m \times \ell \times q}$, $\cB \in \R^{\ell \times r \times q}$ and $\cC \in \R^{r \times n \times q}$, define their product to be the tensor $\cD := \cA * \cB * \cC \in \R^{m \times n \times q}$. Then, it is true that
	\begin{align}
		\widehat{\cD} \equiv \textnormal{\texttt{fft}}\left(\cD,\,[\,\,],\, 3\right) = \chA\, \triangle\, \chB\, \triangle\, \chC
	\end{align}
	where $\chA, \chB$ and $\chC$ are the Fourier domain representations of the tensors $\cA, \cB$ and $\cC$ respectively, with the Fourier transform applied along the third dimension. Further, it holds for a chain of tensor t-product $\cD := \cA_{1}\, *\, \cA_{2}\, *\, \cdots\, *\, \cA_{p}$ that
	$$
		\widehat{\cD} \equiv \textnormal{\texttt{fft}}\left(\cD,\,[\,\,],\, 3\right) = \chA_{1}\, \triangle\, \chA_{2}\, \triangle\, \cdots\, \triangle\, \chA_{p}.
	$$
	\begin{proof}
		We provide a proof only for the triple product case, i.e., $p=3$. The proof for arbitrary $p$ is straightforward. Note that $\cD = \textnormal{\texttt{ifft}}\left( \textnormal{\texttt{fft}}\left( \cD\right) \right) = \textnormal{\texttt{ifft}}\left( \widehat{\cD} \right)$. For the tensor $\cD$, we have
		\begin{align*}
			\cD &= \left(\cA * \cB\right) * \cC,\\
			&= \underbrace{\textnormal{\texttt{ifft}}\left( \textnormal{\texttt{fft}}\left( \cA \right)\, \triangle\, \textnormal{\texttt{fft}}\left(\cB \right) \right)}_{ =: \cZ} *\, \cC, \qquad \text{(see \Cref{rem:fft_multiply})}\\
			&= \cZ *\, \cC,\\
			&= \textnormal{\texttt{ifft}}\left( \textnormal{\texttt{fft}}\left( \textnormal{\texttt{ifft}}\left( \textnormal{\texttt{fft}}\left( \cA \right)\, \triangle\, \textnormal{\texttt{fft}}\left(\cB \right) \right)\, \right)\, \triangle\, \textnormal{\texttt{fft}}\left(\cC \right) \right),\\
			&= \textnormal{\texttt{ifft}}\left( \textnormal{\texttt{fft}}\left( \cA \right)\, \triangle\, \textnormal{\texttt{fft}}\left(\cB \right) \, \triangle\, \textnormal{\texttt{fft}}\left(\cC \right)\, \right).\\
		\end{align*}
		Finally, observing that $\cD = \textnormal{\texttt{ifft}}\left( \widehat{\cD} \right) = \textnormal{\texttt{ifft}}\left( \textnormal{\texttt{fft}}\left( \cA \right)\, \triangle\, \textnormal{\texttt{fft}}\left(\cB \right) \, \triangle\, \textnormal{\texttt{fft}}\left(\cC \right)\, \right)$ completes the proof.
	\end{proof}
\end{lemma}
Next, we state two lemmas on the projection properties of the interpolatory projector $\mathscr{D}$. 
\begin{lemma}
	\label{lem:tqdeim_proj_fourier_domain}
	Suppose $\widehat{\mathscr{D}}$ is the Fourier domain representation of the t-Q-DEIM projector $\mathscr{D} = \cU*(\cP^{\tpose}* \cU)^{-1}*\cP^{\tpose} \in \R^{N \times N \times M}$. Then, each frontal slice of $\widehat{\mathscr{D}}$ is a projector, i.e., 
	$$
			\left(\widehat{\mathscr{D}}^{(i)}\right)^{2} = \widehat{\mathscr{D}}^{(i)} \cdot \widehat{\mathscr{D}}^{(i)} = \widehat{\mathscr{D}}^{(i)}
	$$
	with $\widehat{\mathscr{D}}^{(i)} \in \R^{N \times N}$.
	\begin{proof}
		First, we observe that the assumption $\widehat{\mathscr{D}} \equiv \chU \triangle\, \left(\chP^{\tpose} \triangle\, \chU\right)^{-1} \triangle\, \chP^{\tpose}$ is true by virtue of \Cref{lem:tprod_triple}. Since the product is carried out for each frontal slice, we have that the $i$-th frontal slice is simply $\widehat{\mathscr{D}}^{(i)} := \chU\left(:,\,:,\,i\right)  \left(\chP\left(:,\,:,\,i\right)^{\tpose}  \chU\left(:,\,:,\,i\right)\right)^{-1} \chP\left(:,\,:,\,i\right)^{\tpose}$. Performing, $\widehat{\mathscr{D}}^{(i)} \cdot \widehat{\mathscr{D}}^{(i)}$, the statement can be seen to hold true, since $\left(\chP\left(:,\,:,\,i\right)^{\tpose}  \chU\left(:,\,:,\,i\right)\right)$ is invertible.
	\end{proof}
\end{lemma}

\begin{lemma}
	\label{lem:ident_proj}
	For the t-Q-DEIM projection operator $\mathscr{D} = \cU*(\cP^{\tpose}* \cU)^{-1}*\cP^{\tpose} \in \R^{N \times N \times M}$ and the identity tensor $\cI \in \R^{N \times N \times M}$, the following statement for the t-spectral norm $\norm{\cdot}{}$ is true:
	\begin{align}
		\norm{\cI - \mathscr{D}}{} = \norm{\mathscr{D}}{}.
	\end{align}
	\begin{proof}
		Using the norm equivalence property of the t-spectral norm in the Fourier domain, we have
		$$
		\norm{\cI - \mathscr{D}}{} = \norm{\chI -  \widehat{\mathscr{D}}}{} = \norm{\bcirc{\cI - \mathscr{D}}}{}
		$$
		where $\chI$ and $\widehat{\mathscr{D}}$ are the Fourier domain representations, respectively, of $\cI$ and $\mathscr{D}$. Further, we get
		\begin{align*}
			\norm{\cI - \mathscr{D}}{} &= 
			\left\lVert 
			\begin{pmatrix}
			\chI_{1} - \widehat{\mathscr{D}}_{1}\\
			& \ddots \\
			& & \chI_{M} - \widehat{\mathscr{D}}_{M}
			\end{pmatrix}  
			\right\rVert,\\[1em]
			&= \max \limits_{i = 1, 2, \ldots, M} \norm{\chI_{i} - \widehat{\mathscr{D}}_{i}}{},\\[1em]
			&= \max \limits_{i = 1, 2, \ldots, M} \norm{\widehat{\mathscr{D}}_{i}}{},\\[1em]
			&= \left\lVert 
			\begin{pmatrix}
			\widehat{\mathscr{D}}_{1}\\
			& \ddots \\
			& & \widehat{\mathscr{D}}_{M}
			\end{pmatrix}  
			\right\rVert,\\[1em]
			&= \norm{\bcirc{\mathscr{D}}}{} = \norm{\mathscr{D}}{}.
		\end{align*}
		The third equality above is due to the fact that $\widehat{\mathscr{D}}_{i}$ is a projector (owing to \Cref{lem:tqdeim_proj_fourier_domain}). In the matrix case, $\norm{\mathbf{I} - \mathbf{D}}{2} = \norm{\mathbf{D}}{2}$ for the matrix spectral norm, provided $\mathbf{D}$ is a projector; see~\citep{Szyld06} .
	\end{proof}
\end{lemma} 
Next, we state the main theorem concerning the t-Q-DEIM approximant $\mathfrak{f}_{\text{tq}}$.
\begin{theorem}
	\label{thm:main_theorem}
	Let $\mathfrak{f} \in \R^{N \times 1 \times M}$ be an arbitrary nonlinear, tensor-valued function. Further, given an orthogonal third-order tensor $\cU \in \R^{N \times n \times M}$ and a sampling tensor $\cP \in \R^{N \times n \times M}$ such that $(\cP^{\tpose} * \cU)$ is invertible, let $\mathfrak{f}_{\text{tq}} \in \R^{N \times 1 \times M}$ be the t-Q-DEIM approximation of $\mathfrak{f}$ with $\mathfrak{f}_{\text{tq}} = \cU*(\cP^{\tpose}* \cU)^{-1}*\cP^{\tpose}*\mathfrak{f}$. The error in approximating $\mathfrak{f} $ using $\mathfrak{f}_{\text{tq}}$ can be bounded in the t-spectral norm as follows:
	\begin{equation}
		\norm{\mathfrak{f} - \mathfrak{f}_{\text{tq}}}{} \leq \norm{(\cP^{\tpose} * \cU)^{-1}}{} \cdot \norm{\left(\mathcal{I} - \cU*\cU^{\tpose}\right)* \mathfrak{f}}{}.
	\end{equation}
	\begin{proof}
		We begin by defining the orthogonal projection of $\mathfrak{f}$ in the t-span of $\cU$. We have
		$$
			\mathfrak{f}^{*} := \cU * \cU^{\tpose} * \mathfrak{f}.
		$$
		By the Schmidt-Echart-Young-Mirsky theorem for the t-SVD, $\mathfrak{f}^{*}$ is the optimal approximation of the tensor $\mathfrak{f}$. Next, we define the interpolatory projector $\mathscr{D} := \cU*(\cP^{\tpose}* \cU)^{-1}*\cP^{\tpose}$. Since it is a projector operator, it holds that $\mathscr{D}^{2} = \mathscr{D}$~(\Cref{prop:t-proj}). The t-Q-DEIM approximation then reads
		$$
			\mathfrak{f}_{\text{tq}} = \mathscr{D}*\mathfrak{f}.
		$$ 		
		We observe that
		\begin{equation*}
		\begin{aligned}
			\mathfrak{f} = \left( \mathfrak{f} - \mathfrak{f}^{*} \right) + \mathfrak{f}^{*} = \mathfrak{e}_{f} + \mathfrak{f}^{*}
		\end{aligned}
		\end{equation*}
		where we define the optimal approximation error $\mathfrak{e}_{f} := \mathfrak{f} - \mathfrak{f}^{*}$ with $\mathfrak{e}_{f} \in \R^{N \times 1 \times M}$. Using the definition of the orthogonal projection it holds true that
		\begin{equation}
		\begin{aligned}
			\mathfrak{e}_{f} &= \mathfrak{f} - \mathfrak{f}^{*} = \mathfrak{f} - \cU * \cU^{\tpose} * \mathfrak{f},\\
			&= \left(\cI - \cU*\cU^{\tpose}\right)*\mathfrak{f}.
		\end{aligned}
		\end{equation}
		The t-Q-DEIM approximant $\mathfrak{f}_{\text{tq}}$ can be written as below:
		\begin{equation*}
		\begin{aligned}
			\mathfrak{f}_{\text{tq}} &= \mathscr{D}*\mathfrak{f} = \mathscr{D}*\left(\mathfrak{e}_{f} + \mathfrak{f}^{*}\right),\\
			&= \mathscr{D}*\mathfrak{e}_{f} + \mathscr{D}*\mathfrak{f}^{*}.
		\end{aligned}
		\end{equation*}
		Using the definition of the projector and the optimal approximant $\mathfrak{f}^{*}$
		\begin{equation}
		\label{eq:ftq_expression}
		\begin{aligned}
			\mathfrak{f}_{\text{tq}} &= \mathscr{D}*\mathfrak{e}_{f} + \cU*\left(\cP^{\tpose}*\cU\right)^{-1}*\cP^{\tpose}*(\cU*\cU^{\tpose}*\mathfrak{f}),\\
			&= \mathscr{D}*\mathfrak{e}_{f} + \mathfrak{f}^{*}
		\end{aligned}
		\end{equation}
		where the last equality follows from the fact that $\left(\cP^{\tpose}*\cU\right)$ is invertible.
		For the difference between the true tensor-valued function $\mathfrak{f}$ and its t-Q-DEIM approximation $\mathfrak{f}_{\text{tq}}$, it holds that
		\begin{equation}
		\begin{aligned}
		\mathfrak{f} - \mathfrak{f}_{\text{tq}} &= \mathfrak{e}_{f} - \mathscr{D}*\mathfrak{e}_{f}, \qquad (\text{using}~\cref{eq:ftq_expression})\\
		&= \left(\cI - \mathscr{D}\right) * \mathfrak{e}_{f}.
		\end{aligned}
		\end{equation}
		In the above equation, $\cI \in \R^{N \times N \times M}$ is the t-identity tensor.
		Taking the t-spectral norm on both sides of the above expression and using the sub-multiplicative property~\Cref{lem:t-submultiplicative} yields
		\begin{align*}
			\norm{\mathfrak{f} - \mathfrak{f}_{\text{tq}}}{} &= \norm{\left(\cI - \mathscr{D}\right) * \mathfrak{e}_{f}}{},\\
			&\leq \norm{\cI - \mathscr{D}}{} \cdot \norm{\mathfrak{e}_{f}}{}.
		\end{align*}
		
		Next, we use the fact that $\norm{\cI - \mathscr{D}}{} = \norm{\mathscr{D}}{}$ from \Cref{lem:ident_proj} to write the above inequality as
		\begin{align*}
			\norm{\mathfrak{f} - \mathfrak{f}_{\text{tq}}}{} &\leq \norm{\mathscr{D}}{} \cdot \norm{\mathfrak{e}_{f}}{}.
		\end{align*}		

		Using the definition of the projector $\mathscr{D}$, we obtain
		\begin{align*}
			\norm{\mathfrak{f} - \mathfrak{f}_{\text{tq}}}{} &\leq \norm{\cU*(\cP^{\tpose}* \cU)^{-1}*\cP^{\tpose}}{}\cdot \norm{\mathfrak{e}_{f}}{},\\
			&\leq \norm{\cU}{}\cdot\norm{(\cP^{\tpose}* \cU)^{-1}}{}\cdot\norm{\cP^{\tpose}}{}\cdot \norm{\mathfrak{e}_{f}}{},\\
			&\leq \norm{(\cP^{\tpose}* \cU)^{-1}}{}\cdot \norm{\mathfrak{e}_{f}}{} = \norm{(\cP^{\tpose}* \cU)^{-1}}{}\cdot \norm{\left(\cI - \cU*\cU^{\tpose}\right)*\mathfrak{f}}{}
		\end{align*}		
		where we have utilized the fact the $\norm{\cU}{} = 1$ and $\norm{\cP^{\tpose}}{} = 1$.
	\end{proof}
\end{theorem}

The theorem provides a bound for the approximation error due to t-Q-DEIM in the t-spectral norm. The error bound has two terms in product. The second of these is the optimal-approximation error resulting from an orthogonal projection of the nonlinear tensor $\mathfrak{f}$ on to the t-span of the orthogonal basis $\cU$. The first factor in the bound, viz., $\norm{(\cP^{\tpose}* \cU)^{-1}}{}$ serves as the magnification factor for the optimal error. This is analogous to the standard DEIM error bound, where the quantity $\norm{(\bP^{\tpose} \bU)^{-1}}{2}$ serves as the magnification factor. Intuitively, any procedure to identify the sampling tensor $\cP$ should minimize (or reduce) $\norm{(\cP^{\tpose}* \cU)^{-1}}{}$.

\begin{remark}
	In the standard DEIM and Q-DEIM setting, minimizing the quantity $\norm{(\bP^{\tpose} \bU)^{-1}}{2}$ could be interpreted as being equivalent to reducing the condition number of $\bP^{\tpose} \bU$. To this end, selecting the (most) independent rows of $\bU$ via the pivoted QR decomposition, while sub-optimal, yields good results in practice. In the t-product setting, it is not immediately clear whether the pivoted t-QR based approach~\citep{Haoetal13} we employ to identify the sampling tensor is the optimal choice. It is also not immediately apparent, if the choice of the first frontal slice in the Fourier domain is the best. From our experience, we believe that choosing the frontal slice in the Fourier domain is an appropriate choice as it contains the dominant mode (as a result of taking FFT along the third dimension). 	
	As we illustrate later in the numerical results, this choice yields good approximations in practice for a range of problems. A detailed investigation of other approaches to identify the sampling tensor could be a fruitful task for future research.
\end{remark}

\paragraph{A computable estimate for the orthogonal projection error}
The first quantity in the above bound, viz., $\norm{(\cP^{\tpose}* \cU)^{-1}}{}$ is independent of the argument $\p$ where the function is evaluated. However, the second term ($\norm{\left(\cI - \cU*\cU^{\tpose}\right)*\mathfrak{f}}{}$) depends on $\mathfrak{f}$, and therefore, on $\p$. This leads to repeated evaluations of the second term, on different values of $\p$. Addressing this, we derive a more readily computable estimate, by extending the arguments presented in~\citep{morChaS10} to the tensor t-product setting. 

Let us assume that the tensor $\cF \in \R^{N \times n_{s} \times M}$, whose lateral slices consists of the snapshots matrix of the nonlinear function at different parameters, is representative of the range of $\mathfrak{f}(\p)$ and an $\mathfrak{f}(\p)$ can be approximated as a linear combination of the lateral slices of $\cF$ as
$$
	\mathfrak{f}(\p) = \cF * \mathfrak{h}(\p)
$$
with $\mathfrak{h}(\p) \in \R^{n_{s} \times 1 \times M}$ representing the coefficients of such an approximation. For simplicity, we remove $\p$ from the above equation. Recall from~\cref{eq:snapshot_svd}, the snapshots tensor can be represented in terms of its t-SVD, i.e., $\cF = \cU_{f} * \cS_{f} * \cW_{f}^{\tpose}$, which gives us $\mathfrak{f} = \cU_{f} * \cS_{f} * \cW_{f}^{\tpose}* \mathfrak{h}$. Note that $\cU_{f} \in \R^{N \times N \times M}$, $\cS_{f} \in \R^{N \times n_{s} \times M}$, and $\cW_{f} \in \R^{n_{s} \times n_{s} \times M}$. In~\citep{Kiletal13}, it is shown that the t-SVD can be represented as a summation in the following fashion:
\begin{align*}
	\cF = \cU_{f} * \cS_{f} * \cW_{f}^{\tpose} = \sum\limits_{i=1}^{\min(N, n_{s})} \mathcal{U}_{f}\left(:,\,i,\,:\right) * \cS_{f}\left(i,\,i,\,:\right) * \cW_{f}\left(:,\,i,\,:\right)^{\tpose}.
\end{align*}
Taking this point of view, we express the product $\mathfrak{f} = \cF * \mathfrak{h}$ as
\begin{align*}
	\mathfrak{f} = \sum\limits_{i=1}^{\min(N, n_{s})} \mathcal{U}_{f}\left(:,\,i,\,:\right) * \cS_{f}\left(i,\,i,\,:\right) * \cW_{f}\left(:,\,i,\,:\right)^{\tpose} * \mathfrak{h} = \sum\limits_{i=1}^{\min(N, n_{s})} \mathcal{U}_{f}\left(:,\,i,\,:\right) *\alpha_{i}
\end{align*}
with $\alpha_{i} := \cS_{f}\left(i,\,i,\,:\right) * \cW_{f}\left(:,\,i,\,:\right)^{\tpose} * \mathfrak{h} \in \R^{1 \times 1 \times M}$. Next we observe that
$$
	\cU^{\tpose} * \mathfrak{f} = \sum\limits_{i=1}^{\min(N, n_{s})} \cU\left(:,\,i,\,:\right)^{\tpose} * \mathcal{U}_{f}\left(:,\,i,\,:\right) *\alpha_{i} = \sum\limits_{i=1}^{n} \cI_{n}(:,\,i,\,:) *\alpha_{i}
$$
where $\cI_{n}(:,\,i,\,:) \in \R^{n \times 1 \times M}$ is the $i$-th lateral slice of an identity tensor $\cI_{n} \in \R^{n \times n \times M}$. The second equality follows from the t-orthogonal property~(\Cref{def:t-orth}) of $\cU, \cU_{f}$ (see \Cref{appendix:torth-svd}).
Next, we note that
$$
\cU * \cU^{\tpose} * \mathfrak{f} = \cU * \sum\limits_{i=1}^{n} \cI_{n}(:,\,i,\,:) *\alpha_{i} = \sum\limits_{i=1}^{n} \cU\left(:,\,i,\,:\right) * \alpha_{i}
$$
where the last equality follows from the fact that the t-product of a tensor with the identity retains the tensor. For the difference $\mathfrak{f} - \cU * \cU^{\tpose} * \mathfrak{f}$, it holds
$$
	\mathfrak{f} - \cU * \cU^{\tpose} * \mathfrak{f} = \sum\limits_{i=1}^{\min(N, n_{s})} \cU_{f}\left(:,\,i,\,:\right) * \alpha_{i} - \sum\limits_{i=1}^{n} \cU_{f}\left(:,\,i,\,:\right) * \alpha_{i} = \sum\limits_{i=n+1}^{\min(N, n_{s})}  \cU_{f}\left(:,\,i,\,:\right) * \alpha_{i}.
$$
Expanding $\alpha_{i}$, we have
\begin{align*}
\mathfrak{f} - \cU * \cU^{\tpose} * \mathfrak{f} =  \sum\limits_{i=n+1}^{\min(N, n_{s})} \cU_{f}\left(:,\,i,\,:\right) * \cS_{f}\left(i,\,i,\,:\right) * \cW_{f}\left(:,\,i,\,:\right)^{\tpose} * \mathfrak{h}.
\end{align*}
Defining $\widetilde{\cU} := \cU_{f}\left(:,\,n+1:\min(N, n_{s}),\,:\right)$ to be the tensor containing the last $\min(N, n_{s}) - n$ lateral slices of $\cU_{f}$ and, similarly, defining $\widetilde{\cS} := \cS_{f}\left(n+1:\min(N, n_{s}),\,n+1:\min(N, n_{s}),\,:\right)$ along with $\widetilde{\cW} := \cW_{f}\left(:,\,n+1:\min(N, n_{s}),\,:\right)$, the above summation can be rewritten in a more compact form as shown below:
$$
\mathfrak{f} - \cU * \cU^{\tpose} * \mathfrak{f} = \widetilde{\cU} * \widetilde{\cS} * \widetilde{\cW}^{\tpose} * \mathfrak{h} = \widetilde{\cU} * \widetilde{\cS} * \widetilde{\mathfrak{h}}.
$$
In the above expression, we have defined $\widetilde{\mathfrak{h}} := \widetilde{\cW}^{\tpose} * \mathfrak{h} \in \R^{(\min(N,n_{s})-n) \times 1 \times M}$. 
Taking the t-spectral norm on both sides gives
\begin{align*}
	\norm{\mathfrak{f} - \cU * \cU^{\tpose} * \mathfrak{f}}{} \leq \norm{\widetilde{\cU}}{} \cdot \norm{\widetilde{\cS}}{} \cdot  \norm{\widetilde{\mathfrak{h}}}{} = \norm{\widetilde{\cS}}{} \norm{\widetilde{\mathfrak{h}}}{}
\end{align*}
with the last equality resulting from the observation that $\widetilde{\cU}$ is a t-orthogonal tensor whose t-spectral norm being unity.
A reasonable approximation to the above expression is
$$
\norm{\mathfrak{f} - \cU*\cU^{\tpose}*\mathfrak{f}}{} \lessapprox \norm{\widetilde{\cS}}{}.
$$
Therefore, the spectral norm of the truncated singular tensor, viz., $\widetilde{\cS}$ approximately bounds the orthogonal projection error and serves as a more readily computable estimate to the orthogonal projection error.

\subsection{Discussion concerning the error analysis}
It is evident from the error analysis performed in the previous section that the error in the t-Q-DEIM approximation of a nonlinear tensor-valued data $\mathfrak{f}$ is essentially the norm of the orthogonal projection error multiplied by the amplification factor $\norm{(\cP^{\tpose} * \cU)^{-1}}{}$. Thus, the sampling tensor $\cP$ needs to be chosen such that the norm of this amplification factor is minimized. We further interpret this amplification factor, in light of the definition of the t-spectral norm. We first note that by definition
$$
\norm{(\cP^{\tpose} * \cU)^{-1}}{} =
\left \lVert 
\begin{pmatrix}
\widehat{\cP\cU}^{(1)}\\
& \ddots \\
& & \widehat{\cP\cU}^{(M)}
\end{pmatrix}^{-1}
\right \rVert
$$
with $\widehat{\cP\cU}^{(i)}$ denoting the $i$-th frontal slice of the Fourier domain representation of $(\cP^{\tpose} * \cU) \in \R^{n \times n \times M}$. Specifically, $\widehat{\cP\cU} := \texttt{fft}((\cP^{\tpose} * \cU),\, [\,\,],\,3)$.
Since the matrix is block diagonal, the inverse is applied block-wise, resulting in
\begin{align*}
\norm{(\cP^{\tpose} * \cU)^{-1}}{} 
&= 
\left \lVert 
\begin{pmatrix}
\left(\widehat{\cP\cU}^{(1)}\right)^{-1}\\
& \ddots \\
& & \left(\widehat{\cP\cU}^{(M)}\right)^{-1}
\end{pmatrix}
\right \rVert,\\
&= \max \limits_{i = 1, 2, \ldots, M} \left\lVert \left(\widehat{\cP\cU}^{(i)}\right)^{-1}\right\rVert.
\end{align*}
The last equality is an outcome of the application of the t-spectral norm to each block. Moreover, note that each $\left(\widehat{\cP\cU}^{(i)}\right)^{-1} \in \C^{n \times n}$ is a matrix. For any matrix $\mathbf{A}$, it is known that $\norm{\bA^{-1}}{2} = 1 / \sigma_{\text{min}}(\bA)$. Based on this, we obtain
\begin{align*}
\norm{(\cP^{\tpose} * \cU)^{-1}}{}
&= \max \limits_{i = 1, 2, \ldots, M} \dfrac{1}{ \sigma_{\text{min}}\left(\widehat{\cP\cU}^{(i)}\right)}.
\end{align*}
Finally, we note that $\left(\widehat{\cP\cU}^{(i)}\right) = \left(\widehat{\cP}\widehat{\cU}^{(i)}\right)$ since $\widehat{\cP}$ has all its frontal slices to be the same, being a projector. Then the goal of the sampling algorithm is to find a $\widehat{\cP}$ in order to minimize the t-spectral norm of the amplification factor, via the equivalent objective function
\begin{align}
 \min\,\, \norm{(\cP^{\tpose} * \cU)^{-1}}{} &=  \min \limits_{\chP} \max \limits_{i = 1, 2, \ldots, M} \dfrac{1}{\sigma_{\text{min}}\left(\widehat{\cP}\widehat{\cU}^{(i)}\right)}.
\end{align}
The above min-max reformulation of the problem of identifying the `best' index locations chosen by the tensor $\cP$ seems not amenable to an analytical solution. The brute-force approach would cost $M$ (or $M/2$ for real-valued tensors) separate QR decompositions for the frontal slices followed by the same number of SVDs of the sampled $n \times n$ matrix corresponding to each frontal slice. In contrast, our method of choosing the first frontal slice alone is cheaper and also yields good approximation, at least for the five examples considered in this work.

\paragraph{A-priori bound on the amplification factor}
An a-priori upper bound for the magnification factor $\norm{(\bP^{\tpose} \bU)^{-1}}{2}$ was provided in the original DEIM approach~\citep{morChaS10}. This bound was further sharpened in~\citep{morDrmG16}. The Q-DEIM upper bound for the magnification factor is based on the dimensions of the basis matrix $\bU \in \C^{N \times n}$ and reads
\begin{align}
	\label{eq:qdeim_apriori_bound}
	\norm{(\bP^{\tpose} \bU)^{-1}}{2} \leq \dfrac{\sqrt{N - n + 1}}{\sigma_{\text{min}}(\bU)} \dfrac{\sqrt{4^{n} + 6n - 1}}{3}
\end{align}
with $\sigma_{\text{min}}(\bU)$ representing the smallest singular value of $\bU$. The main argument in obtaining the above bound employed by~\citep{morDrmG16} is that the pivoting operation on the $\bU^{T}$ matrix results in a diagonally dominant leading block of the (upper triangular) $\mathbf{R}$ matrix resulting from the QR decomposition of $\bU^{T}$. Standard arguments on norms of triangular matrices are then used to derive the bound.

To develop an a priori bound for $\norm{(\cP^{\tpose}* \cU)^{-1}}{}$, we note that
$$
	\norm{(\cP^{\tpose}* \cU)^{-1}}{} = \norm{(\chP^{\tpose} \triangle\, \chU)^{-1}}{},
$$
that is, the t-spectral norm of the magnification factor can be equivalently computed based on its Fourier domain counterparts. Recall that the triangle $\triangle$ represents the frontal slice-wise matrix multiplication operation. We further note that each frontal slice of $\chU \in \R^{N \times n \times M}$ is Hermitian. Since the sampling tensor is computed based on the QR-decomposition of the first frontal slice of $\chU$, the Q-DEIM a-priori bound holds, i.e.,
$$
	\norm{\left(\chP(:,\,:,\,1)^{\tpose} \, \chU(:,\,:,\,1)\right)^{-1}}{} \leq \dfrac{\sqrt{N - n + 1}}{\sigma_{\text{min}}(\chU(:,\,:,\,1))} \dfrac{\sqrt{4^{n} + 6n - 1}}{3}.
$$
Owing to the frontal slices in the Fourier domain being Hermitian, this entails $\sigma_{\text{min}}(\chU(:,\,:,\,j)) = 1$. However, in case of t-Q-DEIM, it is not straightforward to provide an a-priori bound for the full amplification factor $\norm{(\cP^{T} * \cU)^{-1}}{}$ as in the case of the Q-DEIM. This is because the sampling tensor is chosen based on the first frontal slice. As every frontal slice is independent, the Q-DEIM a-priori bound~\cref{eq:qdeim_apriori_bound} valid for the first frontal slice is no longer applicable to the remaining ones. An a-priori error estimator that seemingly works in practice is
\begin{align}
\label{eq:apriori_bound}
\norm{\left(\chP^{\tpose} \, \chU\right)^{-1}}{} \lessapprox M \cdot \dfrac{\sqrt{N - n + 1}}{1} \dfrac{\sqrt{4^{n} + 6n - 1}}{3}.
\end{align}
\begin{remark}
We believe this estimate can be tightened further (or even made a rigorous bound) by following the arguments in the Q-DEIM bound and arriving at an equivalent expression through the diagonal dominance property of the $\mathcal{R}$ tensor from the t-pQR decomposition. We detail the computation of the t-pQR decomposition in \Cref{appendix:tpqr}. We hypothesize that a rigorous a-priori upper bound for the amplification factor would leverage the diagonal dominance structure of the $\mathcal{R}$ tensor resulting from the t-pQR decomposition.
\end{remark}

\subsection{Oversampling for better approximation and robustness to noise}
Real-world data is often noisy. Sparse approximation approaches are often susceptible to noise in the data, leading to poor approximation quality. Prior work has shown that oversampling the function at a limited amount of extra points often delivers robustness to noise~\citep{PehDG20}. Such an approach can be extended to tensor data as well. Note that, since the t-Q-DEIM basis $\cU$ is of dimension $n$, we choose $n$ sampling locations, this is the interpolation regime. By oversampling the horizontal slices of $\cU$ (as done for the basis matrix $\bU$ in previous work), additional sampling points can be obtained. By doing so, we are no longer in the interpolation regime. The t-Q-DEIM approximation with oversampled data is then a best least-squares approximation. The coefficient $\mathfrak{c}$ is obtained as
\begin{align}
\mathfrak{c} = \left(\cP^{\tpose} * \cU \right)^{\dagger} * \cP^{\tpose} * \mathfrak{f}(\p)
\end{align} 
which gives us the oversampled t-Q-DEIM approximation
$$
\mathfrak{f}_{\text{tq-OS}} = \cU * \left(\cP^{\tpose} * \cU \right)^{\dagger} * \cP^{\tpose} * \mathfrak{f}(\p)
$$
with $\cU \in \R^{N \times n \times M}$ and $\cP \in \R^{N \times n_{O} \times M}$, where $n_{O} > n$. This requires additional considerations concerning the computation of the pseudo-inverse in the tensor t-product and also requires extending some classical singular value bounds to the case of the t-SVD. As this is beyond the immediate scope of the present work, it will be addressed in a forthcoming publication.

\subsection{Training and inference costs}
\label{subsec:training_and_inference_costs}
We briefly discuss the computational costs for the training and inference stages of both t-Q-DEIM and Q-DEIM. For simplicity, we use the same notation for the reduced dimension $n$ for both t-Q-DEIM and Q-DEIM. In practice, they can be different. Often, $n_{\text{t-Q-DEIM}} < n_{\text{Q-DEIM}}$. For ease of notation we define the following computational costs for operations involving third-order tensors:
\begin{itemize}
	\item The cost of FFT / Inverse FFT of a tensor $\cA \in \R^{n_{1} \times n_{2} \times n_{3}}$ is
	$$
		\mathcal{C}_{\text{fft}}(n_{1}, n_{2}, n_{3}) := \mathcal{O}(n_{1} n_{2} n_{3} \log(n_{3})),
	$$
	\item The cost of t-product of two tensors $\cA \in \R^{n_{1} \times n_{2} \times n_{4}}$ and $\cB \in \R^{n_{2} \times n_{3} \times n_{4}}$ is denoted $\mathcal{C}_{\text{t-pro}d}(n_{1}, n_{2}, n_{3}, n_{4})$ and amounts to
	$$
	 \mathcal{C}_{\text{t-prod}}(n_{1}, n_{2}, n_{3}, n_{4}) := \mathcal{C}_{\text{fft}}(n_{1}, n_{2}, n_{4}) + \mathcal{C}_{\text{fft}}(n_{2}, n_{3}, n_{4}) + \mathcal{O}(n_{1} n_{2} n_{3} n_{4}).
	$$
\end{itemize}

\subsubsection{Training cost}
The dominant cost in the training stage for both t-Q-DEIM and Q-DEIM is the SVD. 

\paragraph{t-Q-DEIM}For the t-Q-DEIM, data tensor needs to be first transformed into the Fourier domain. For the data tensor $\cFtrain \in \R^{m \times N_{\text{train}} \times q}$, this incurs cost $\mathcal{C}_{\text{fft}}(m,N_{\text{train}},q)$. In the Fourier domain, a separate matrix SVD is performed for every frontal slice, amounting to a total of $q$ separate SVDs on matrices of dimension $\R^{m \times N_{\text{train}}}$. If the data tensor $\cFtrain$ is real-valued , this cost can be reduced in half by taking advantage of conjugate symmetry in the Fourier domain and performing only half the amount of SVDs. This results in cost $0.5 \cdot q \mathcal{O}(\max(m, N_{\text{train}} ), \min(m, N_{\text{train}})^{2})$. Following this, the resulting SVD factors are truncated to dimension $n$. Performing an inverse FFT on the truncated factors has cost $\mathcal{C}_{\text{fft}}(m,n,q)$.
 The pivoted QR factorization is done on the first frontal slice of the basis tensor $\cU \in \R^{m \times n \times q}$ in the Fourier domain. Typically, for the examples we considered, $m > n$. Therefore, in the worst case, the QR decomposition has complexity $\mathcal{O}(n^{3})$. Computing the matrix $\mathcal{D}$ through a tensor-tensor product of tensors $\cU \in \R^{m \times n \times q}$ and $(\cU(\mathbf{p},\,:,\,:))^{-1} \in \R^{n \times n \times q}$ has cost $\mathcal{C}_{\text{t-prod}}(m, n, n, q)$
 for the tensor-tensor product and $0.5 q \mathcal{O}(n^{3}) + 2 	\mathcal{C}_{\text{fft}}(m, n, q)$ for the tensor inverse.
Summing up, the total offline cost for the t-Q-DEIM adds up to
\begin{align*}
	\mathcal{C}_{\text{t-Q-DEIM}} &:= \phantom{ }\mathcal{C}_{\text{fft}}(m,N_{\text{train}},q) \qquad\qquad (\texttt{FFT})\\ &+ 0.5 q \mathcal{O}(\max(m, N_{\text{train}} ), \min(m, N_{\text{train}})^{2}) \qquad (\texttt{SVD})\\ &+ \mathcal{C}_{\text{fft}}(m,n,q) \qquad (\texttt{IFFT})\\&+ \mathcal{O}(n^{3}) \qquad (\texttt{QR})\\ &+\mathcal{C}_{\text{t-prod}}(m, n, n, q) + 0.5 q \mathcal{O}(n^{3}) + 2 \mathcal{C}_{\text{fft}}(m, n, q)  \qquad (\texttt{compute}\,  \mathcal{D})
\end{align*}

\paragraph{Q-DEIM}For the Q-DEIM algorithm, the standard matrix SVD can be achieved in at most $\mathcal{O}(\max(m, q N_{\text{train}} ), \min(m, q N_{\text{train}})^{2})$. Depending on the application, $m > q N_{\text{train}}$ or $m < q N_{\text{train}}$. Following this, the QR decomposition of the basis $\mathbf{U} \in \R^{m \times n}$ has worst-case complexity $\mathcal{O}(n^{3})$. Finally, evaluating $\mathbf{D}$ incurs cost scaling as $\mathcal{O}(mn^{2})$ for the matrix-matrix product and $\mathcal{O}(n^{3})$ for the dense matrix inversion.

Summing up, the total offline cost for the Q-DEIM adds up to
\begin{align*}
\mathcal{C}_{\text{Q-DEIM}} &:= \mathcal{O}(\max(m, q N_{\text{train}} ), \min(m, q N_{\text{train}})^{2}) \qquad (\texttt{SVD})\\
&+ \mathcal{O}(n^{3}) \qquad (\texttt{QR})\\
&+ \mathcal{O}(mn^{2}) + \mathcal{O}(n^{3})  \qquad (\texttt{compute}\,  \mathbf{D})
\end{align*}

\subsubsection{Inference cost}
Computing the cost of performing inference is straightforward. We assume that samples of the target function is available at $n$ spatial locations. For the t-Q-DEIM method we have data samples $\mathfrak{f}_{\text{sampled}} \in \R^{n \times 1 \times q}$. Computing $\mathfrak{f}_{\text{tq}}$ from this involves a tensor-tensor product which can be achieved with cost $\mathcal{C}_{\text{t-prod}}(m, n, 1, q)$. Note that, in case $q$ dimension denotes time, we achieve inference at all time locations in one-shot. 

For the Q-DEIM, data is assumed also to available at $n$ spatial locations, i.e., $\bff_{\text{sampled}} \in \R^{n}$. Computing the Q-DEIM approximation has computational complexity $\mathcal{O}(m n)$. To evaluate $q$ different time instances, the cost is $\mathcal{O}(m q n)$. Here, it is critical to note that $n$ for t-Q-DEIM will be different from $n$ for Q-DEIM. For a desired level of approximation accuracy, we often get $n_{\text{t-Q-DEIM}} < n_{\text{Q-DEIM}}$. Thus, the inference cost of t-Q-DEIM can still be comparable to Q-DEIM. 

\section{Numerical results}
\label{sec:numerics}
We illustrate the proposed t-Q-DEIM approach on five numerical examples. These are:
\begin{enumerate}
	\item the Burgers' equation \hfill \textsf{(1-parameter, 1D, nonlinear)}
	\item the FitzHugh-Nagumo equations \hfill \textsf{(2-parameter, 1D, nonlinear)}
	\item the thermal cookie problem \hfill \textsf{(4-parameter, 2D, linear)}
	\item the Navier-Stokes equations (flow past a square cylinder) \hfill \textsf{(1-parameter, 2D, nonlinear)}
	\item the brain interface dataset \hfill \textsf{(Experimental data)}
\end{enumerate}
The examples are carefully chosen to illustrate the robust approximation properties of the t-Q-DEIM method for systems characterized by different dimensions, number of parameters, and\\(non-)linearity. We not only consider spatio-temporal data from parametrized dynamical systems, but also test on an experimental dataset arising in biology. For each example we consider, the data is assumed to be arranged in a data tensor $\cF \in \R^{m \times \ell \times q}$. We perform interpolation along the lateral slice, viz., the second dimension of the third-order tensor. The data tensor is further divided into training and testing datasets, $\cF_{\text{train}} \in \R^{m \times N_{\text{train}} \times q}$ and $\cF_{\text{test}} \in \R^{m \times N_{\text{test}} \times q}$, respectively, with $\ell = N_{\text{train}} + N_{\text{test}}$. The t-Q-DEIM algorithm is applied to the training data $\cF_{\text{train}}$ and the performance of the approximation is tested on making predictions/forecast on $\cF_{\text{test}}$. For comparison, we also test the performance of the Q-DEIM approach on both these datasets.

\subsection{Code availability}
The codes to reproduce the numerical experiments will be made available upon publication.

\subsection{Note on the computational environment}
All numerical results were obtained on a desktop computer running Ubuntu 20.04, installed with a $12$-th generation \INTEL \textsf{i5} processor, $32$GB of RAM. 
The code is written in Python, using the Spyder IDE. For ready reference, the computational timings for each example to be discussed is provided in~\Cref{tab:timings}. We compare runtimes of both the training and inference stages for t-Q-DEIM~(\Cref{alg:tqdeim}) and Q-DEIM~(\Cref{alg:qdeim}). For every example, the values provided in~\Cref{tab:timings} are the mean over $k$ independent runs. To clearly see the dominant role of the SVD computations, the time taken for the SVD is provided separately.

\begin{table}[b!]
	\centering
	\def\arraystretch{2}\tabcolsep=10pt
	\small
\begin{tabular}{|c|cccccc|}
	\hline
	\multirow{3}{*}{Example}                 & \multicolumn{3}{c|}{t-Q-DEIM runtime}                                                                      & \multicolumn{3}{c|}{Q-DEIM runtime}                                                                        \\ \cline{2-7} 
	& \multicolumn{2}{c|}{Offline}                        & \multicolumn{1}{c|}{\multirow{2}{*}{Online}} & \multicolumn{2}{c|}{Offline}                        & \multicolumn{1}{l|}{\multirow{2}{*}{Online}} \\ \cline{2-3} \cline{5-6}
	& \multicolumn{1}{l}{SVD} & \multicolumn{1}{c|}{Rest} & \multicolumn{1}{c|}{}                        & \multicolumn{1}{l}{SVD} & \multicolumn{1}{c|}{Rest} & \multicolumn{1}{l|}{}                        \\ \hline
	Burgers' ($n=8$, $k=10$)                 & $1.30$                  & $2.55$e-2                 & $2.47$e-2                                    & $2.35$                  & $3$e-4                    & $3$e-3                                       \\ \hline
	FitzHugh-Nagumo ($n=21$, $k=10$)         & $0.96$                  & $5.64$e-1                 & $3.56$e-2                                    & $2.91$                  & $6.38$e-4                 & $9.83$e-4                                    \\ \hline
	Thermal cookie ($n=35$, $k=5$)           & $13.03$                 & $1.30$                    & $3.22$e-2                                    & $95.10$                 & $6.55$e-3                 & $2.21$e-3                                    \\ \hline
	Navier-Stokes ($n=20$, $k=5$)            & $7.59$                  & $2.19$                    & $1.2$e-1                                     & $155.80$                & $3.45$e-3                 & $7.8$e-3                                     \\ \hline
	Brain interface dataset ($n=13$, $k=10$) & $8.26$e-2               & $2.9$e-3                  & $9.46$e-5                                    & $1.82$e-1               & $6.53$e-5                 & $2.03$e-6                                    \\ \hline
\end{tabular}
	\caption{Offline and online computation times for the t-Q-DEIM and Q-DEIM approaches. The values in the brackets next to each example indicate the reduced dimension ($n$) and number of independent runs ($k$) over which the reported average runtimes were computed. All values are in seconds.}
	\label{tab:timings}
\end{table}

\subsection{Burgers' equation}
\begin{figure}[t!]
	\centering
	\includegraphics[scale=0.7]{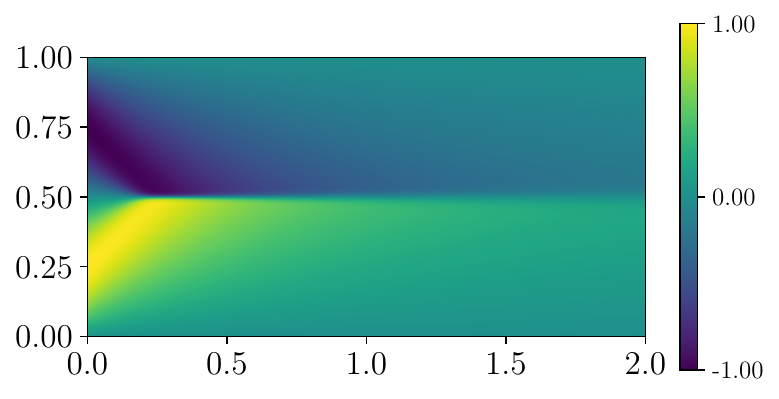}
	\caption{Solution to the viscous Burgers' equation at viscosity $\mu = 0.004$.}
	\label{fig:burgers_fom}
\end{figure}

\begin{figure}[h!]%
	\centering
	\includegraphics[scale=0.6]{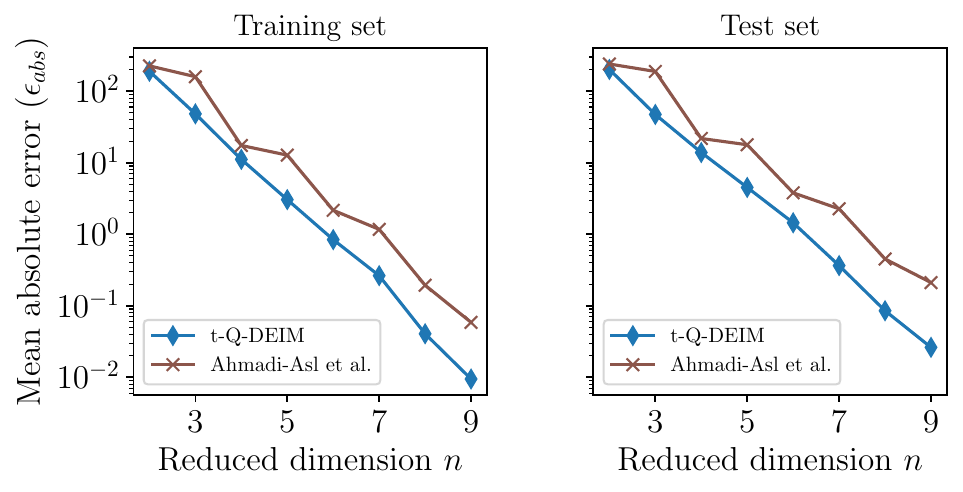}
	\caption{Burgers' equation: influence of the sampling strategy on the approximation errors; sampling approach from~\Cref{alg:tpqr} and the method proposed in~\citep{Asletal24} are compared.}
	\label{fig:burgers_compare_point_selections}
\end{figure}

\paragraph{Model description}
We consider the viscous Burgers' equation defined as
\begin{align}
	\dfrac{\partial}{\partial t} w(x, t) &= w(x, t) \dfrac{\partial}{\partial x} w(x, t) + \mu \dfrac{\partial^{2}}{\partial x^{2}} w(x, t)
\end{align}
with boundary conditions $w(x=0, t;\mu) = w(x=1, t;\mu) = 0$ and initial condition $w(x, t~=~0;\mu) = \sin(x)$. Here, $w(x, t)$ is the state variable of interest. The spatial variable $x \in [0,\,1]$ and the time $t \in [0,\,2]$. The above PDE was discretized in space with the second-order finite difference method (dimension $N = 1000$) and further discretized in time using a first-order implicit-explicit (IMEX) method into $N_{t} = 1000$ time steps. In~\Cref{fig:burgers_fom}, we plot the space-time solution of the Burgers' equation at $\mu = 0.004$. To obtain the data tensor $\cF$, the discretized system was simulated at $N_{\text{s}} = 50$ log-uniformly spaced parameter samples of the viscosity $\mu \in [0.004,\,0.04]$. The $50$ samples were then randomly divided into $N_{\text{train}} = 25$ and $N_{\text{test}} = 25$ samples, respectively. The training data tensor $\cFtrain \in \R^{1000 \times 25 \times 1000}$ and the test data tensor $\cFtest \in \R^{1000 \times 25 \times 1000}$. In effect, each lateral slice of the training (testing) data tensor consists of $N_{t}=1000$ solution snapshots of the viscous Burgers' equation at a given parameter $\mu$, at different time instances.
\paragraph{Comparison to a standard DEIM like approach for index selection}
Recently, in the context of CUR decomposition of tensors, the work~\citep{Asletal24} introduces a standard-DEIM like slice selection algorithm. Since it is a closely related approach to ours, we first make a comparison of the respective performances of our proposed t-Q-DEIM approach and the method from~\citep{Asletal24}. We perform a sensitivity analysis-type experiment, plotting the approximation error in the t-spectral norm as a function of the reduced dimension $n$. For the Burgers' equation dataset, we consider $n \in \{2, 3, 4, \ldots, 9\}$. The analysis is carried out for both t-Q-DEIM (\Cref{alg:tqdeim}) and the method from~\citep{Asletal24}. We first observe that the proposed t-Q-DEIM approach takes $20.16$ seconds for the entire analysis where $8$ separate approximations were obtained. Roughly, this translates to $2.52$ seconds for generating a single approximation. The same analysis using the other method takes $20.82$ seconds, equivalent to $2.60$ seconds for producing each approximation. Thus, we see that the proposed t-pQR based point selection strategy is computationally slightly faster, even for this small example. We believe that the slightly higher run time of the method in~\citep{Asletal24} is likely due to the sequential or iterative choice of the sampling indices, where multiple norm computations need to be performed.  In~\Cref{fig:burgers_compare_point_selections}, the error convergence plots for the training and testing data are shown. 
For each approximation $\widetilde{\mathfrak{f}}(\p)$, we obtain the mean value of the true error over the training and test sets, viz. 
$$
\epsilon_{\text{abs}} := \dfrac{1}{N_{\text{samples}}} \sum\limits_{i=1}^{N_{\text{samples}}} \| \mathfrak{f}(\p) - \widetilde{\mathfrak{f}}(\p) \|
$$
with $\p$ taken from the training set or the testing set and $N_{\text{samples}}$ either $N_{\text{train}}$ or $N_{\text{test}}$. 
We can infer immediately that the proposed t-Q-DEIM sampling strategy yields better approximation quality consistently. 

\begin{figure}[h!]
	\centering
	\includegraphics[scale=0.6]{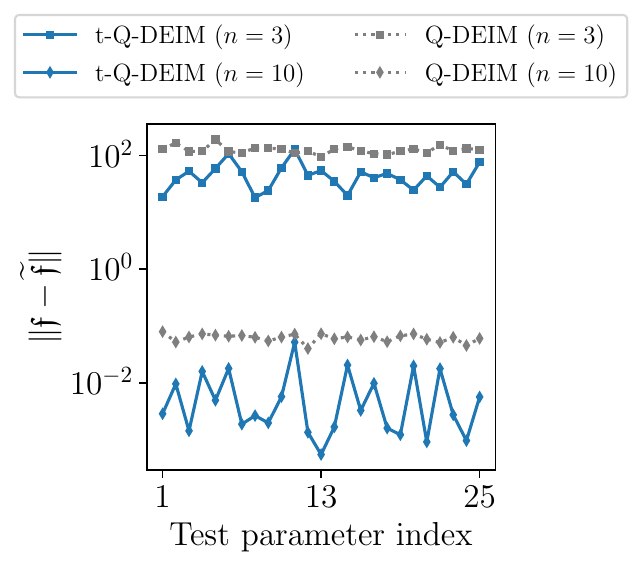}
	\caption{Burgers' equation: comparison of the performance of t-Q-DEIM and Q-DEIM on the test parameter set for $n = 3$, $n = 10$, respectively.}
	\label{fig:burgers_test_set_error_compare}
\end{figure}

\begin{figure}[t]%
	\centering
	\begin{subfigure}[t]{0.45\linewidth}
		\centering
		\includegraphics[scale=0.75]{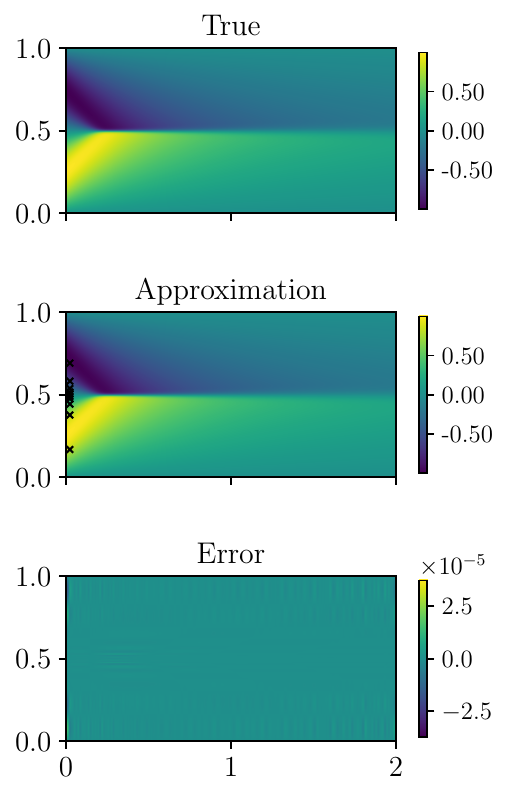}
		\caption[a]{$\p = 0.00439$}
	\end{subfigure}\hfill
	\begin{subfigure}[t]{0.45\linewidth}
		\centering
		\includegraphics[scale=0.75]{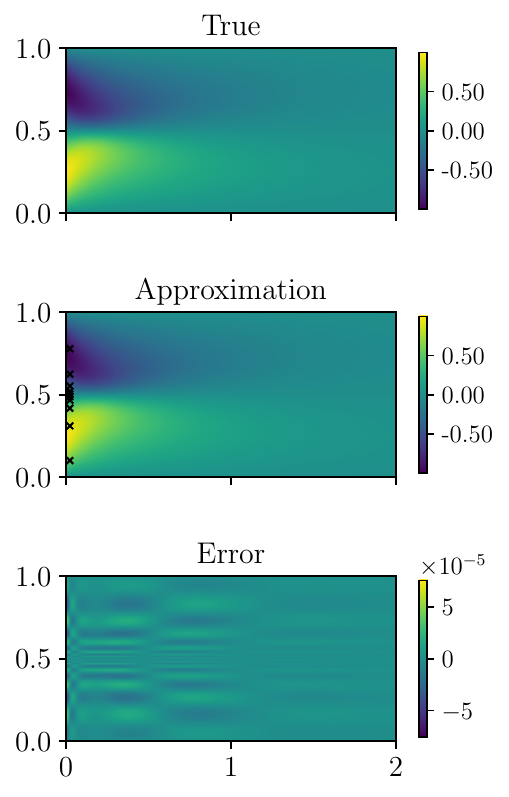}
		\caption[b]{$\p = 0.04$}
	\end{subfigure}
	\caption{t-Q-DEIM approximation of the Burgers' equation; each figure plots, from top to bottom, the true solution, the t-Q-DEIM approximation, and the pointwise errors. The black cross marks denote the location where the data is sampled.}
	\label{fig:burgers_tqdeim_test_soln}
\end{figure}

\paragraph{Performance of t-Q-DEIM compared to Q-DEIM}
We apply both the proposed t-QDEIM (\Cref{alg:tqdeim}) and Q-DEIM (\Cref{alg:qdeim}) to the Burgers' equation dataset. Both approaches are trained using $\cFtrain$. To evaluate their performance, we illustrate their respective errors (in the t-spectral norm) for two values of $n$, viz. $n = 3$ and $n = 10$ in \Cref{fig:burgers_test_set_error_compare}. The x-axis displays the indices of the parameter samples in the test set, while the y-axis has the error of the respective approximation in the (t-)spectral norm. t-Q-DEIM displays a better performance over Q-DEIM, yielding up to one order of magnitude improvement.  In~\Cref{fig:burgers_tqdeim_test_soln}, the true space-time solution to the Burgers' equation at $\p = 0.00439$ and $\p = 0.04$ (both taken from the test set) and the corresponding t-Q-DEIM solution for $n=10$ are plotted. It can be observed that the t-Q-DEIM approximation is uniformly good, with errors of order $10^{-5}$. Moreover, both low and high-viscosity solutions are captured effectively. The crosses on the plots of the approximation indicate the locations of the t-Q-DEIM samples. Notice that the samples are clustered around the area where the shock in the solution develops. With reference to the computational cost, the training stage of both t-Q-DEIM and Q-DEIM are dominated by the cost of the SVD. For the t-Q-DEIM, the t-SVD requires $1.3$ seconds, while the cost of performing the matrix SVD on the vectorized dataset is nearly double for the Q-DEIM, taking $2.35$ seconds. However, the inference cost for the Q-DEIM approach is more competitive.

\subsection{FitzHugh-Nagumo equations}
\label{num:fhn}
The next example we consider is the 1-D FitzHugh-Nagumo equations. The system models the neuronal activations in response to external stimuli and finds application in brain modeling and also in cardiac electrophysiology. It is known to exhibit a characteristic limit cycle behaviour for certain parameter choices. 
\paragraph{Model description}
The FitzHugh-Nagumo equations, defined in a 1-D spatial domain $x \in \Omega := [0, 1]$, are a coupled system of nonlinear time-dependent equations given by
\begin{equation}
\label{eq:fhn-pde}
\begin{aligned}
	\epsilon \dfrac{\partial}{\partial t} w_{1}(x, t) &= \epsilon^{2} \dfrac{\partial^{2}}{\partial x^{2}} + g_{1}(w_{1}(x, t)) - w_{2}(x, t) + c,\\
	\dfrac{\partial}{\partial t} w_{2}(x, t) &= b w_{1}(x, t) - \gamma w_{2}(x, t) + c
\end{aligned}
\end{equation}
and with the boundary conditions
$$
	\dfrac{\partial}{\partial x} w_{1}(0, t) = -I_{\text{ext}}(t), \qquad \dfrac{\partial}{\partial x} w_{1}(L, t) = 0,
$$ 
and initial condition
$$
	 w_{1}(x, 0) = 0.001, \qquad  w_{2}(x, 0) = 0.001.
$$
The variable $w_{1}(x, t) \in \R^{}$ is the neuron electric potential and  $w_{2}(x, t) \in \R^{}$ is the recovery rate of the potential. $x \in \Omega$ denotes the spatial variable and the time variable $t \in [0, 5]$. The constitutive relation $g_{1}(w_{1}(x, t)) := w_{1} (w_{1} - 0.1) (1 - w_{1})$ represents the nonlinear term. The external stimulus to the system is given by an exponentially decaying input of the form $I_{\text{ext}}(t) := 50000 t^{3} e^{-15 t}$. Potentially, the variables $\epsilon, b, c, \gamma$ can be parameters of the FitzHugh-Nagumo system. We fix $b = 0.5$ and $\gamma = 2$ and treat the system as having two free parameters $\p := (\epsilon, c) \in R := [0.01, 0.04] \times [0.025, 0.075]$. The PDE is spatially discretized with a second-order finite difference scheme with $512$ nodes for each of the two couple equations, yielding a discretized system of dimension $N = 1024$. Similar to the Burgers' equation, time discretization is carried out with a first-order IMEX scheme into $501$ time nodes with $\Delta t = 0.01$. 
\begin{figure}[t]%
	\centering
	\begin{subfigure}[t]{0.56\linewidth}
		\centering
		\includegraphics[scale=0.7]{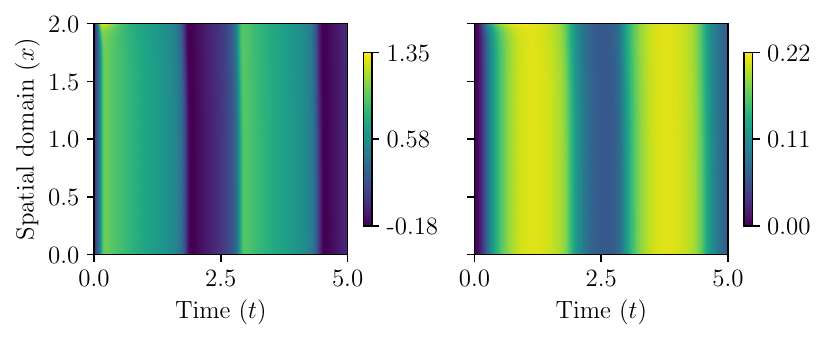}
		\caption[a]{Snapshots of variables $w_{1}, w_{2}$ }
	\end{subfigure}\hfill
	\begin{subfigure}[t]{0.41\linewidth}
		\centering
		\includegraphics[scale=0.505]{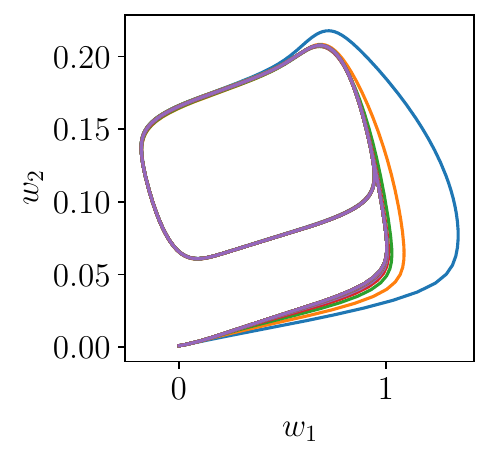}
		\caption[b]{Limit cycle behaviour}
	\end{subfigure}
	\caption{Snapshots corresponding to the variables $w_{1}, w_{2}$ in~\cref{eq:fhn-pde} and the limit cycle exhibited by the FitzHugh-Nagumo equations at parameter $\p := (\epsilon, c) = (0.022, 0.075)$}
	\label{fig:fhn_fom}
\end{figure}
\Cref{fig:fhn_fom} illustrates the solution to the FitzHugh-Nagumo system at a particular value of the parameters $(\epsilon, c)$. As noted before, for some parameter combinations, the system exhibits a limit cycle behaviour.

\begin{figure}[h!]%
	\centering
	\includegraphics[scale=0.7]{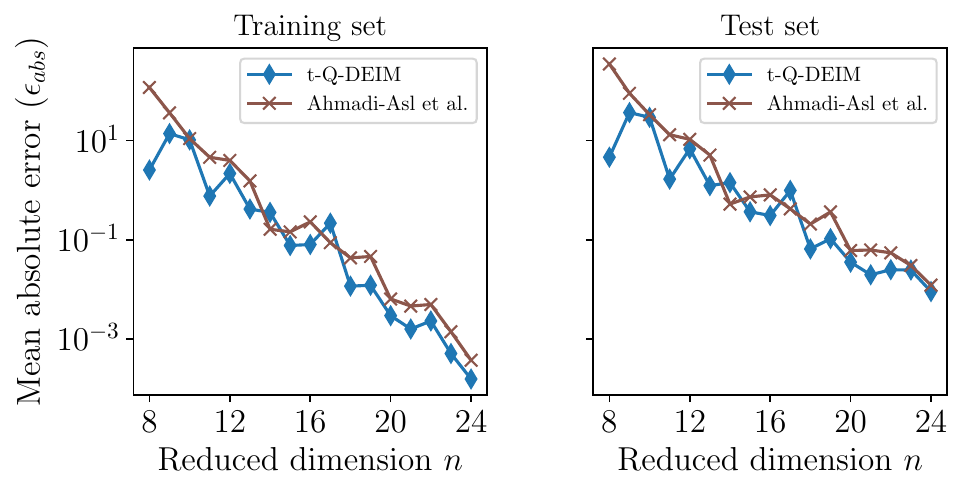}
	\caption{FitzHugh-Nagumo equations: influence of the sampling strategy on the approximation errors; sampling approach from~\Cref{alg:tpqr} and the method proposed in~\citep{Asletal24} are compared.}
	\label{fig:fhn_compare_point_selections}
\end{figure}

The discretized FitzHugh-Nagumo system is solved at $N_{s} = 36$ parameters sampled uniformly from the parameter space $R$. We achieve this by creating a 2D Cartesian grid and picking $6$ equally-spaced samples for each of the two parameters. The data tensor $\cF \in \R^{1024 \times 36 \times 501}$ is then divided into training and test data as $\cFtrain \in \R^{1024 \times 26 \times 501}$, containing $N_{\text{train}} = 26$ lateral slices, and $\cFtest \in \R^{1024 \times 10 \times 501}$, containing $N_{\text{test}} = 10$ lateral slices, with each lateral slice (in the test and training data) having the solution snapshot vector (of dimension $N=1024$) at $N_{t} = 501$ time steps.
\paragraph{Comparison to a standard DEIM like approach for index selection}
As done for the previous Burgers' equation example, we compare the t-pQR based sampling index selection strategy~\Cref{alg:tpqr} with the approach proposed in~\cite{Asletal24}. We carry out a sensitivity analysis as before for $n \in \{8, 9, 10, \ldots, 22, 23, 24\}$, a total of $17$ different approximations. The t-Q-DEIM approach took $34.28$ seconds for this analysis, which is about $2.02$ seconds per approximation. The same implementation for the other sampling strategy takes $52.60$ seconds or $3.09$ seconds per approximation. Evidently, the sampling strategy of~\citep{Asletal24} is more expensive. Next, we plot the error convergence plots to check which of the two strategies yields better performance. We see that in~\Cref{fig:fhn_compare_point_selections}, the t-pQR approach results in a better approximation error, both over the training and test set. Indeed, for the training set, the approach from~\citep{Asletal24} gives a better error at two values of $n$, $n = 14, 17$. Additionally, this strategy also provides slightly better performance on the test set at these same values of $n$. Notwithstanding this, taking together the lesser computational cost of the t-pQR strategy and its overall better performance, we can safely conclude that the t-pQR based sampling strategy outperforms the other method.

\begin{figure}[b!]
	\centering
	\includegraphics[scale=0.6]{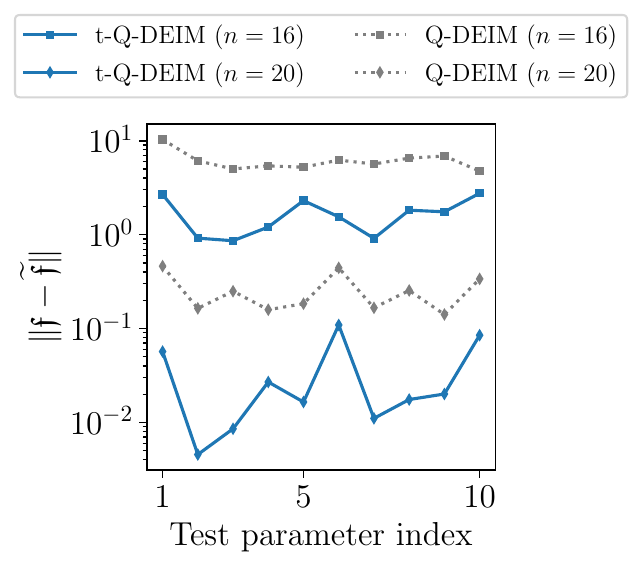}
	\caption{FitzHugh-Nagumo equations: comparison of the performance of t-Q-DEIM and Q-DEIM on the test set for $n = 11$, $n = 20$.}
	\label{fig:fhn_test_set_error_compare}
\end{figure}

\begin{figure}[h!]
	\centering
	\includegraphics[scale=0.6]{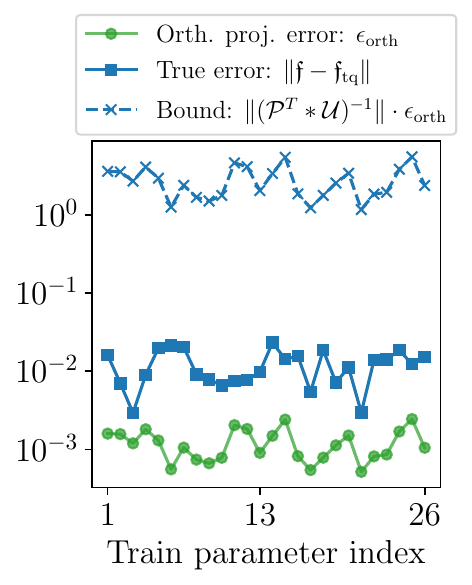}
	\caption{FitzHugh-Nagumo equation: true error, best orthogonal projection error, and the proposed error bound at $n = 19$.}
	\label{fig:fhn_verify_error_bound}
\end{figure}

\begin{figure}[h!]
	\centering
\begin{subfigure}[h!]{1\linewidth}
	\centering
	\includegraphics[scale=0.6]{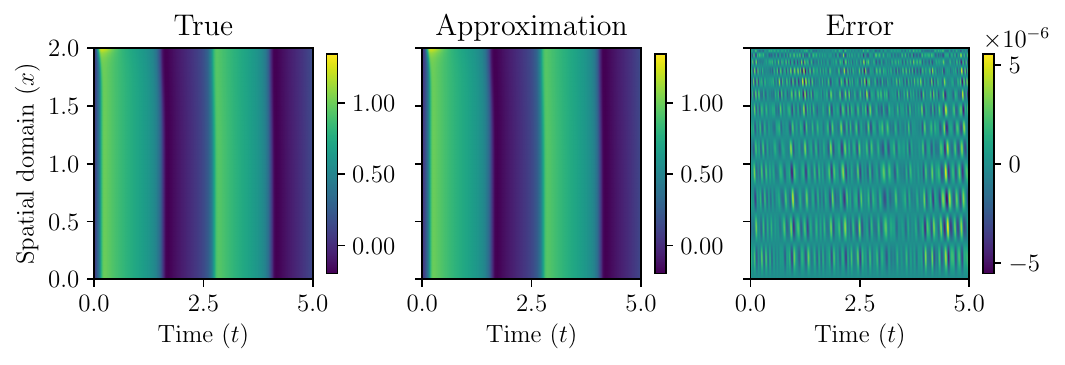}
	\caption[a]{t-Q-DEIM}
\end{subfigure}\hfill
\begin{subfigure}[h!]{1\linewidth}
	\centering
	\includegraphics[scale=0.6]{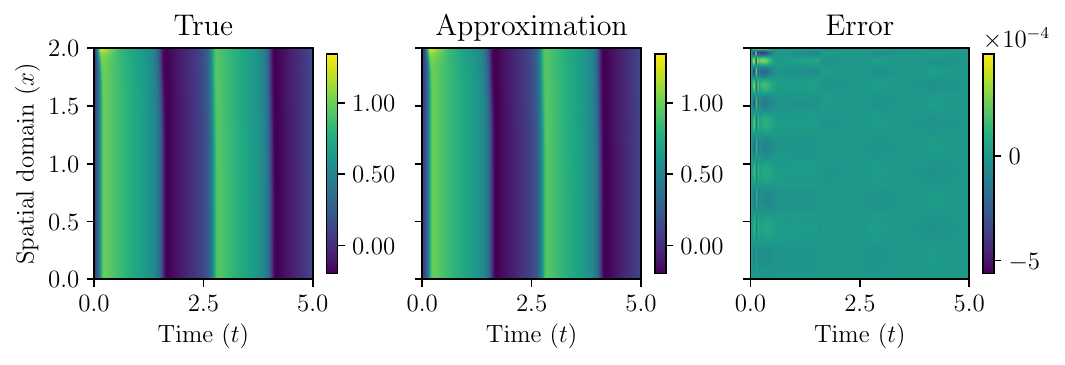}
	\caption[b]{Q-DEIM}
\end{subfigure}
\caption{t-Q-DEIM and Q-DEIM approximation of the FitzHugh-Nagumo equations at $\p := (0.022. 0.065)$; each figure plots, from left to right, the true solution, the approximation, and the pointwise errors.}
\label{fig:fhn_tqdeim_qdeim_test_soln}
\end{figure}

\paragraph{Performance of t-Q-DEIM}
We use the training data $\cFtrain$ to learn the corresponding t-Q-DEIM quantities $\mathcal{U}, \mathcal{P}$ needed for the approximation. As evident from \Cref{fig:fhn_test_set_error_compare}, the proposed t-Q-DEIM is able to achieve considerably better approximation on the test set when compared to Q-DEIM. Next, we verify the validity of the upper bound proposed in \Cref{thm:main_theorem}. Recall that, the error incurred by t-Q-DEIM is essentially the orthogonal projection error on to the subspace of $\cU$ multiplied by a factor. The factor was the t-spectral norm of $(\cP^{T} * \cU)^{-1}$. We plot the three quantities, viz., the true error $\| \mathfrak{f} - \mathfrak{f}_{\text{tq}} \|$, the orthogonal projection error $\| \left(\cI - \cU*\cU^{\tpose}\right) * \mathfrak{f} \|$, and the upper bound, $\| (\cP^{T} * \cU)^{-1} \| \cdot \| \left(\cI - \cU*\cU^{\tpose}\right) * \mathfrak{f} \| $ in \Cref{fig:fhn_verify_error_bound}. The t-Q-DEIM approximation error is about an order of magnitude larger than the best approximation error achievable by an orthogonal projector. More importantly, the bound is uniformly applicable over the entire training set, thus validating~\Cref{thm:main_theorem}.
Next, we plot the t-QDEIM and Q-DEIM approximations of the FitzHugh-Nagumo solution at the test parameter $\p := (0.022. 0.065)$ in \Cref{fig:fhn_tqdeim_qdeim_test_soln}. The t-Q-DEIM yields a better approximation, with an error that is two orders of magnitude less, when compared to Q-DEIM. The offline training cost for the FitzHugh-Nagumo example is dominated by the cost of performing the SVD (t-SVD)~\Cref{tab:timings}. In this aspect, the t-SVD is cheaper ($0.96$ seconds) as compared to the matrix SVD needed for Q-DEIM ($2.91$ seconds), which is almost $3$ times the former. Nevertheless, the Q-DEIM approximation is more competitive in the inference stage, with an improved performance of slightly more than one order of magnitude. This boils down to the cheaper matrix vector product for Q-DEIM over the tensor-tensor product in the case of t-Q-DEIM.

\subsection{Thermal cookie example}
\begin{figure}[t!]
	\centering
	\includegraphics[scale=0.65]{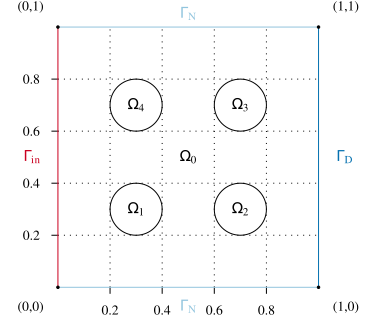}
	\caption{Computational domain of the thermal cookie problem~\citep{morwiki_thermalblock}}
	\label{fig:cookie}
\end{figure}

The thermal cookie example solves the linear heat equation in a square domain, with four inner circular patches, each exhibiting different heat conductivities; see~\Cref{fig:cookie}.
\begin{figure}[tp]%
	\centering
	\begin{subfigure}[t]{0.45\linewidth}
		\centering
		\includegraphics[scale=0.55]{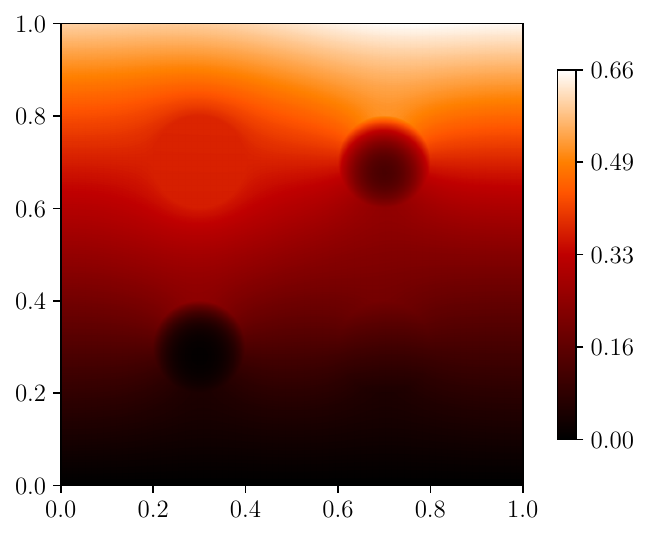}
		\caption[a]{$t = 0.3$s, $\p = (21.5443, 0.0046, 0.1, 0.01)$}
	\end{subfigure}\hfill
	\begin{subfigure}[t]{0.45\linewidth}
		\centering
		\includegraphics[scale=0.55]{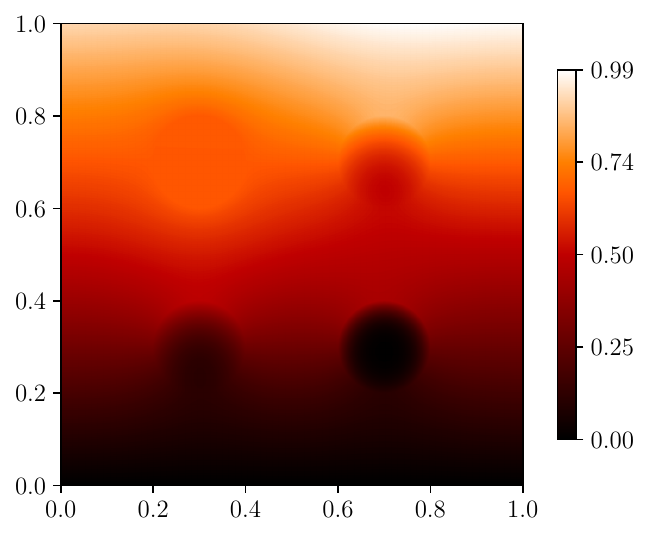}
		\caption[b]{$t = 0.8$s, $\p = (100.0, 0.0046, 0.001, 0.01)$}
	\end{subfigure}
	\caption{Snapshot of the thermal cookie problem at two different parameter and time instances.}
	\label{fig:thermal_fom}
\end{figure}
\paragraph{Model description}
The computational domain of interest is $\Omega := (0,1)^{2}$ which is divided into five subdomains denoted $\{\Omega_{i}\}_{i=0}^{4}$ with  $\Omega = \Omega_{0} \cup \Omega_{1} \cup \Omega_{2} \cup \Omega_{3} \cup	 \Omega_{4}$. For details on the geometry, we refer to the work~\citep{morwiki_thermalblock}.
The governing PDE is given by
\begin{align}
	\dfrac{\partial }{\partial t} \bw(\bx, t; \p) &= -\nabla \cdot (-\p(x) \nabla \bw(\bx, t; \p) )
\end{align}
and it is defined within the domain $\Omega$. The PDE takes in an input applied on the left boundary $\Gamma_{\text{in}} := \{0\} \times (0, 1)$ and
$$
	\p(x) \nabla \bw(\bx, t; \p) \cdot n(\bx) = u(t)
$$
for $t \in [0, T]$, $T = 1$ and $x \in \Gamma_{\text{in}}$. We further have Neumann and Dirichlet conditions on the remaining boundaries. The Neumann boundary is defined as $\Gamma_{\text{N}} := (0,1) \times {0, 1}$ while the Dirichlet boundary is $\Gamma_{\text{D}} := \{1\} \times (0,1)$. The Neumann boundary condition is 
$$
		\p(x) \nabla \bw(\bx, t; \p) \cdot n(\bx) = 0, \qquad x \in \Gamma_{\text{N}}
$$ and the Dirichlet boundary condition is
$$
	\bw(\bx, t; \p) = 0, \qquad x \in \Gamma_{\text{D}}.	
$$
The system parameter is $\p := (\mu_{1}, \mu_{2}, \mu_{3}, \mu_{4}) \in \R^{4}$ with $\mu_{i}$ referring to the heat conductivities on the respective domains $\Omega_{i}$, $i=1,2,3,4$. Furthermore, $\mu_{1} \in [1.0,\,100.0]$, $\mu_{2}, \mu_{3} \in [0.001,\,0.1]$ and $\mu_{4} \in [0.0001,\,0.01]$.
The heat conductivity in $\Omega_{0}$ is uniformly set to unity.
The PDE is discretized with $N = 7488$ nodes using finite element method (using FEniCS 2019.1). The resulting ordinary differential equations are further discretized in time $t \in [0, 1]$ into $N_{t} = 101$ time steps. \Cref{fig:thermal_fom} plots two sample solutions to the thermal cookie problem, one at $\p = (21.5443, 0.0046, 0.1, 0.01)$ and time $t = 0.3$ seconds and the other at $\p = (100.0, 0.0046, 0.001, 0.01)$ and time $t = 0.8$ seconds. Each solution exhibits considerable variation spatially, both in its magnitude and distribution. This makes the thermal problem a particularly challenging one, despite its linear nature.

For this problem, we construct the data tensor $\cF$ by choosing $N_{\text{s}} = 4^{4} = 256$ samples by forming a Cartesian grid consisting of $4$ log-spaced samples for each parameter $\mu_{i}$. This data tensor is then divided further into training and testing matrices (in the ratio $80\,:\,20$); the training matrix has dimension $\cFtrain \in \R^{7488 \times N_{\text{train}} \times 101}$ with $N_{\text{train}} = 200$ (the parameters chosen randomly). Moreover, this results in the testing matrix $\cFtest \in \R^{7488 \times N_{\text{test}} \times 101}$ with $N_{\text{test}} = 56$. 

\begin{figure}[t!]
	\centering
	\includegraphics[scale=0.6]{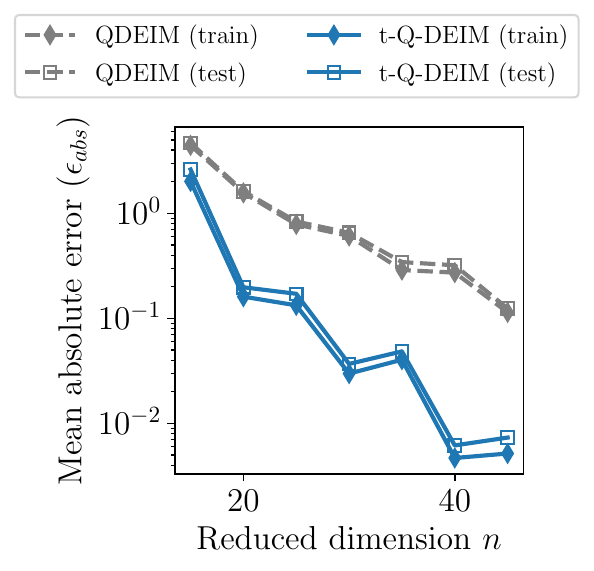}
	\caption{Thermal cookie example: Mean estimated errors over training (testing) set vs. $n$.}
	\label{fig:thermal_sensitivty_analysis}
\end{figure}

\begin{figure}[t!]
	\centering
	\begin{subfigure}[h!]{1\linewidth}
		\centering
		\includegraphics[scale=0.7]{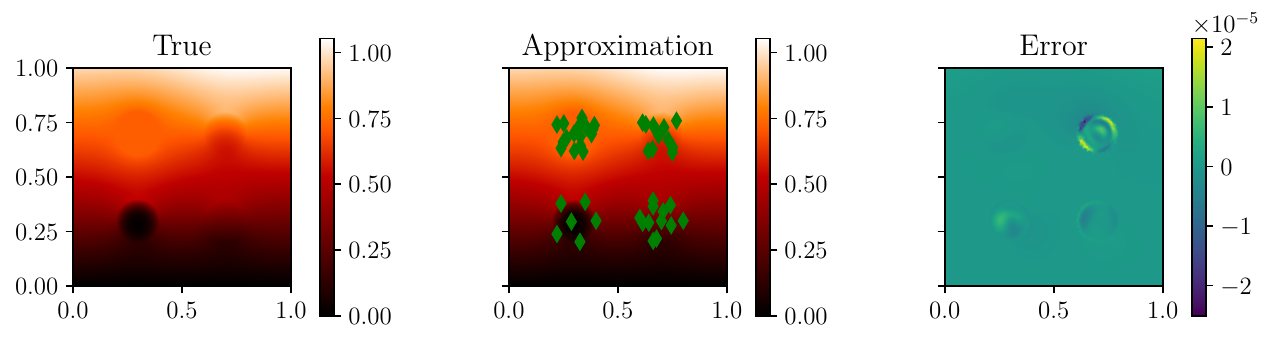}
		\caption[a]{t-Q-DEIM}
	\end{subfigure}\hfil
	\begin{subfigure}[h!]{1\linewidth}
		\centering
		\includegraphics[scale=0.7]{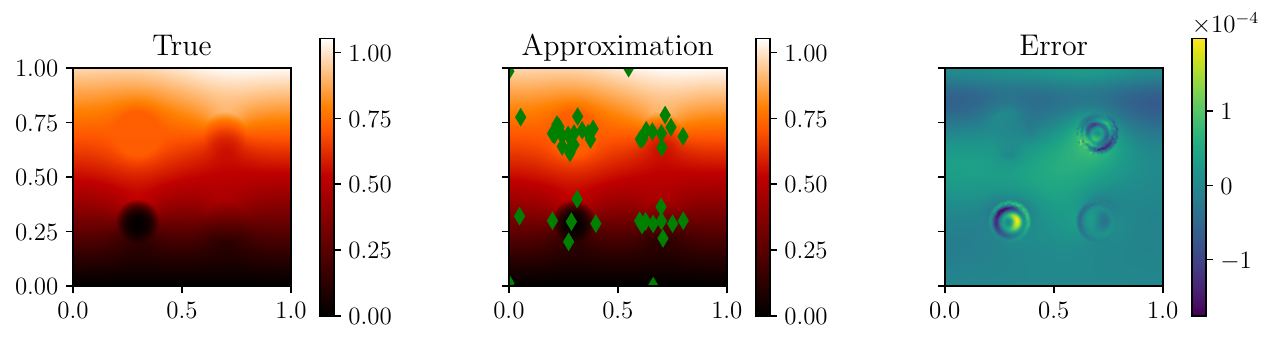}
		\caption[b]{Q-DEIM}
	\end{subfigure}
	\caption{t-Q-DEIM and Q-DEIM approximation of the thermal cookie problem at $\p := (100, 0.01, 0.0215, 0.01)$ at the final time $T = 1.0$ seconds; each figure plots, from left to right, the true solution, the approximation, and the pointwise errors.}
	\label{fig:thermal_tqdeim_qdeim_test_soln}
\end{figure}

\begin{figure}[h!]%
	\centering
	\begin{subfigure}{0.45\linewidth}
		\centering
		\includegraphics[scale=0.6]{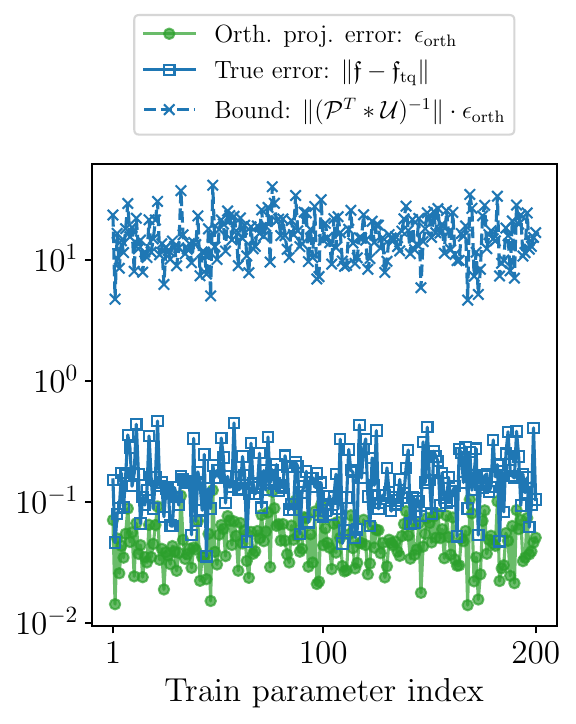}
		\caption[a]{$n = 20$}
	\end{subfigure}\hfill
	\begin{subfigure}{0.45\linewidth}
		\centering
		\includegraphics[scale=0.6]{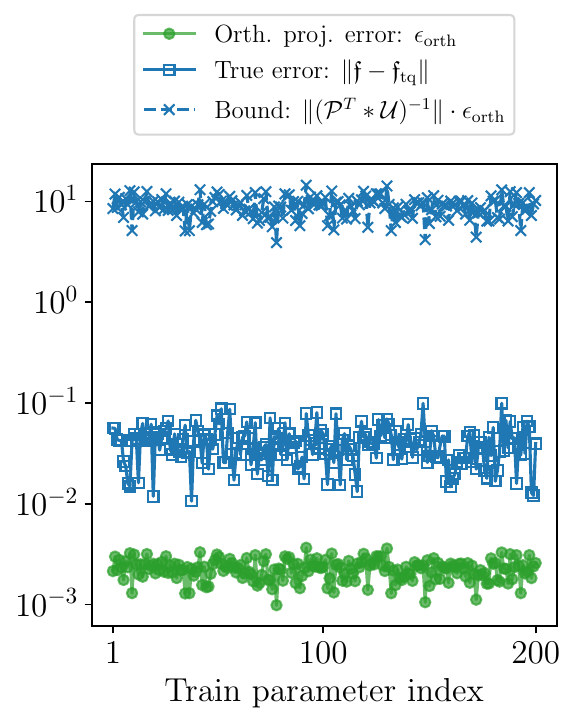}
		\caption[b]{$n = 35$}
	\end{subfigure}
	\caption{Thermal cookie example: true error, best orthogonal projection error, and the proposed error bound.}
	\label{fig:thermal_verify_error_bound}
\end{figure}

\paragraph{Performance of t-Q-DEIM and Q-DEIM} For this example, we begin by performing a sensitivity analysis of the t-Q-DEIM and Q-DEIM approximations. We make use of~\Cref{alg:tqdeim,alg:qdeim}, respectively, on the training data tensor $\cFtrain$ to obtain the approximations at different values of the reduced dimension $n$. The results are plotted in~\Cref{fig:thermal_sensitivty_analysis}. The proposed t-Q-DEIM approach emerges as the clear winner, with higher approximation accuracy for the training and testing parameters, for this multi-parameter system.
Having established the supremacy of the t-Q-DEIM for this example, we move on to verifying the proposed error bound. In \Cref{fig:thermal_verify_error_bound}, the error bound is illustrated for reduced dimensions $n = 20$ and $n = 35$. In both instances, we numerically verify that the upper bound property holds uniformly over all training parameters. It can also be seen that the t-Q-DEIM error is close to the best approximation error when $n = 20$ and is around one order of magnitude larger for $n = 35$. Finally, for this example, we show the Q-DEIM and t-Q-DEIM approximations of the solution at a test parameter $\p := (100, 0.01, 0.0215, 0.01)$ at the final time $T = 1.0$ seconds. As visualized in~\Cref{fig:thermal_tqdeim_qdeim_test_soln}, the proposed method of t-Q-DEIM possesses a better approximation, a fact already borne out through~\Cref{fig:thermal_sensitivty_analysis}. It is interesting to note the locations where t-Q-DEIM and Q-DEIM enforce a sampling point (visualized as green diamonds in the figure). In the case of the former, most of the samples are concentrated around the circular patches inside the domain. This is reasonable to expect as those regions display variations of the thermal conductivity and t-Q-DEIM is able to correctly identify those spatial locations (and the points closely surrounding them) as being critical. In the case of Q-DEIM, while many samples are indeed clustered around the circular patches, a few points are also present close to the boundary of the domain. 

The computational benefits of t-Q-DEIM over Q-DEIM are more starkly visible for this example; see~\Cref{tab:timings}. While the t-SVD needs just $13.03$ seconds, the cost of the matrix SVD for Q-DEIM is $95.1$ seconds, which is nearly seven times worse. Yet, it is indeed observed that, as for the previous examples, the inference costs for the Q-DEIM are close to one order of magnitude cheaper for Q-DEIM.

\subsection{Navier-Stokes equation}
\begin{figure}[t!]
	\centering
	\includegraphics[scale=0.5]{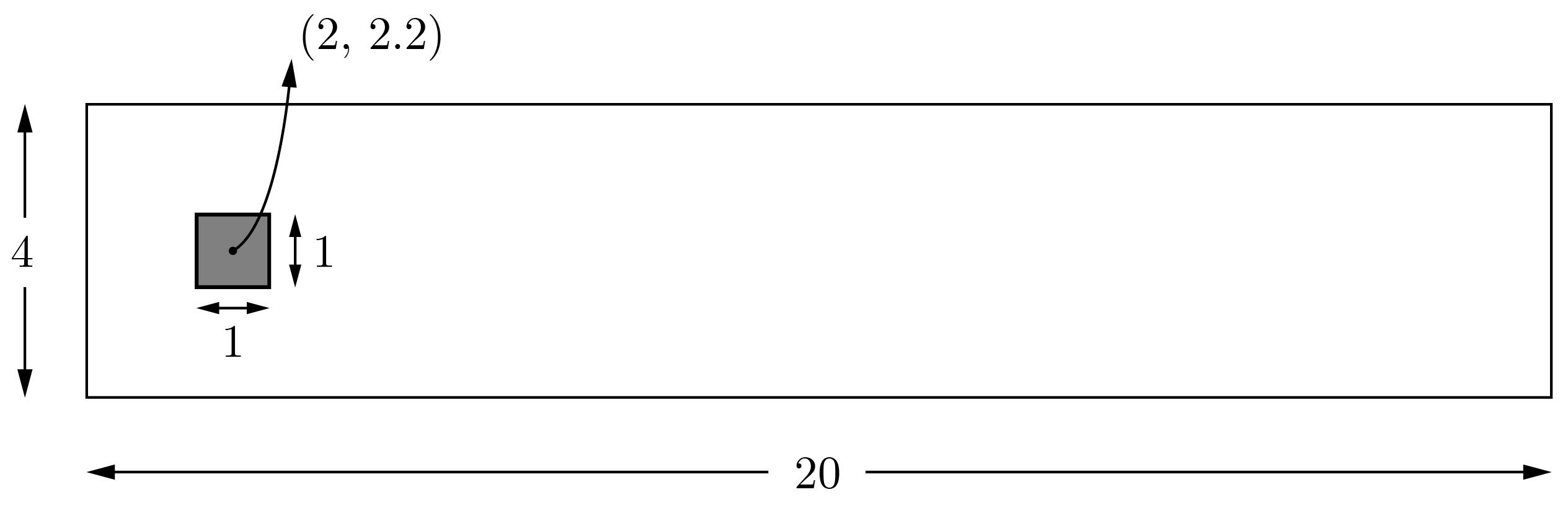}
	\caption{2-D geometry for solving Navier-Stokes equation}
	\label{fig:geom_ns}
\end{figure}
\begin{figure}[b!]
	\centering
	\includegraphics[scale=0.65]{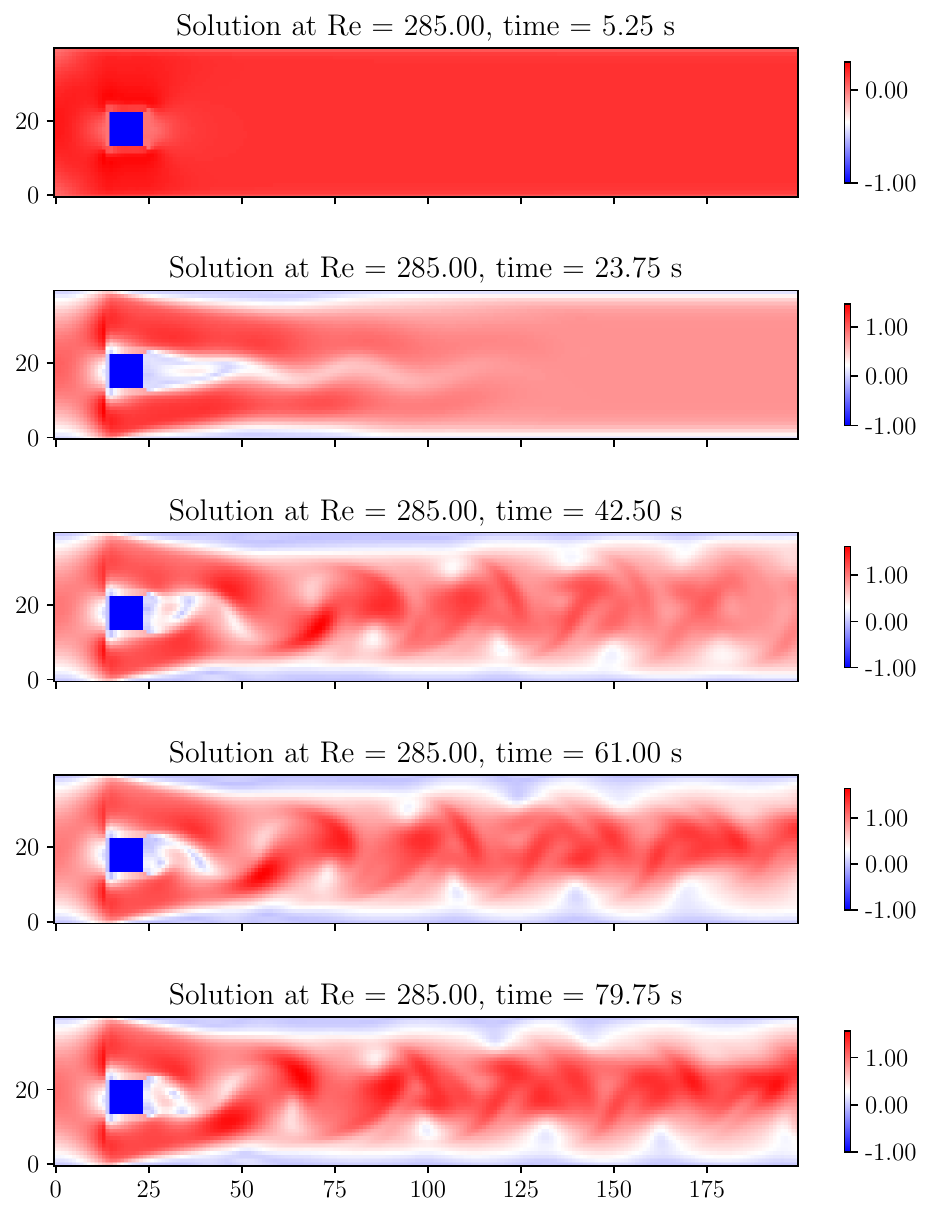}
	\caption{Snapshots of the 2-D Navier-Stokes equation --- flow past a square cylinder at $Re = 285$.}
	\label{fig:ns_fom}
\end{figure}

We next consider the Navier-Stokes equations solved inside the 2-D geometry shown in~\Cref{fig:geom_ns}. 
\paragraph{Model description}
The domain consists of a rectangular channel of dimension $(20 \times 4)$ with a $(1 \times 1)$ solid square obstacle. The PDE is solved using a monolithic finite volume method and an incremental pressure correction scheme~\citep{fvm_repo}. The spatial discretization is done with a uniform 2-D grid, having grid size $0.1$ and a step size of $\Delta t= 0.025$ is used for the time discretization with the time $t \in [0, 80]$ seconds. This results in a discretized system of dimension $N = 40\cdot200 = 8000$. The parameter of interest is the Reynolds' number $\p := \text{Re} \in [105,\,400]$. The solution to the Navier-Stokes system at Reynolds numbers $Re = 285$ is plotted for different time instances $t \in \{5.25,\, 23.75,\, 42.5,\, 61.00,\, 79.75\}$ seconds in~\Cref{fig:ns_fom}. 

\begin{figure}[t!]
	\centering
	\includegraphics[scale=0.6]{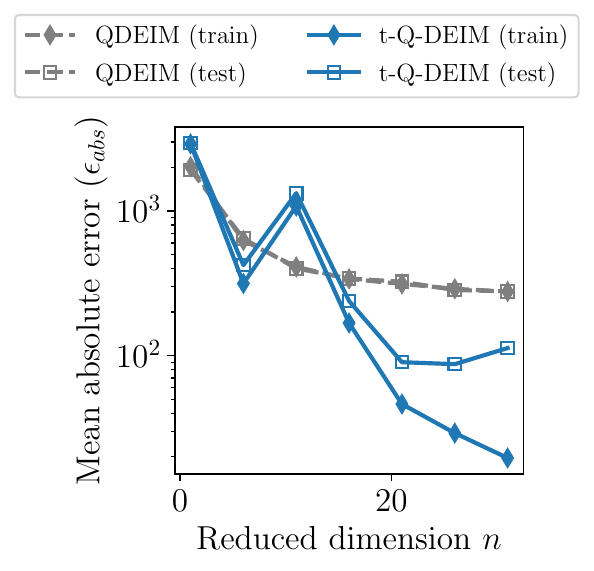}
	\caption{Navier-Stokes equation: Mean estimated errors over training (testing) set vs. $n$.}
	\label{fig:ns_sensitivty_analysis_abserr}
\end{figure}

\begin{figure}[h!]%
	\centering
	\begin{subfigure}{0.45\linewidth}
		\centering
		\includegraphics[scale=0.6]{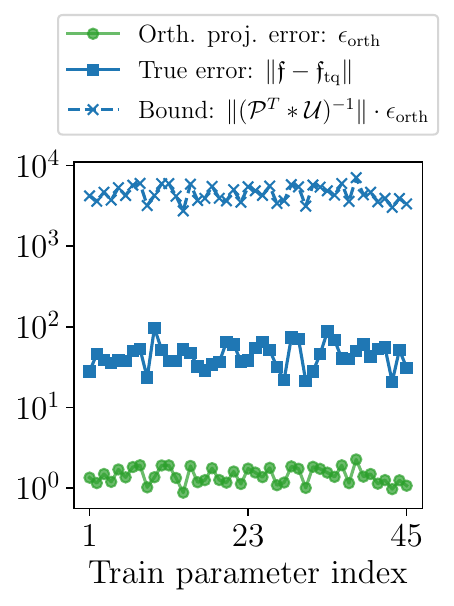}
		\caption[a]{$n = 21$}
	\end{subfigure}\hfill
	\begin{subfigure}{0.45\linewidth}
		\centering
		\includegraphics[scale=0.6]{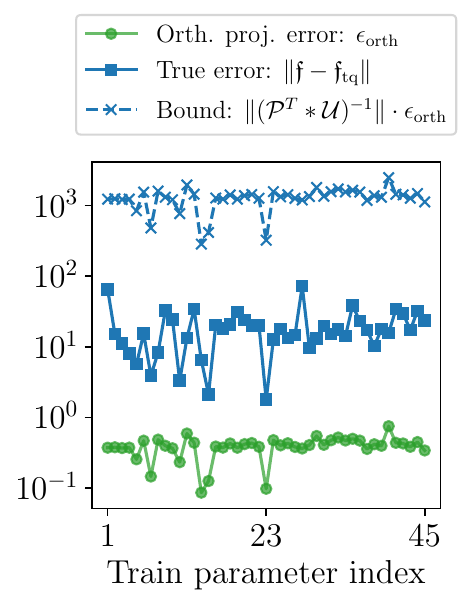}
		\caption[b]{$n = 31$}
	\end{subfigure}
	\caption{Navier-Stokes equation: true error, best orthogonal projection error, and the proposed error bound.}
	\label{fig:ns_verify_error_bound}
\end{figure}

\begin{figure}[t!]
	\centering
	\begin{subfigure}[h!]{0.45\linewidth}
		\centering
		\includegraphics[scale=0.63]{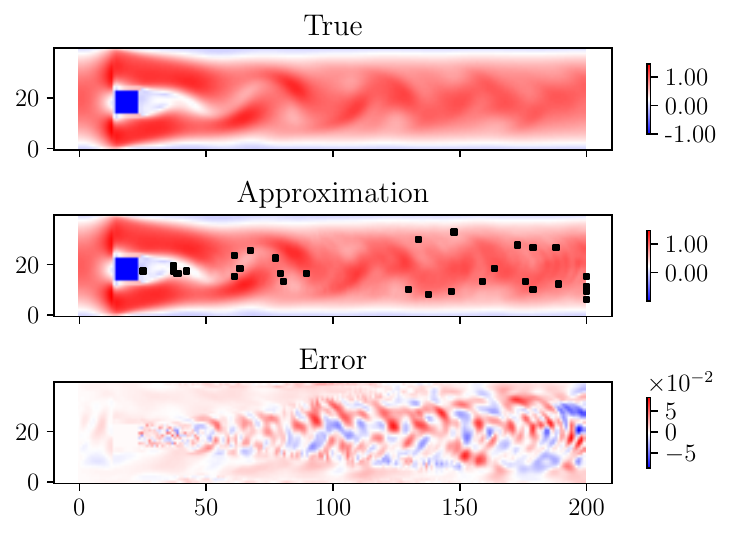}
		\caption[a]{t-Q-DEIM}
	\end{subfigure}\hfill
	\begin{subfigure}[h!]{0.45\linewidth}
		\centering
		\includegraphics[scale=0.63]{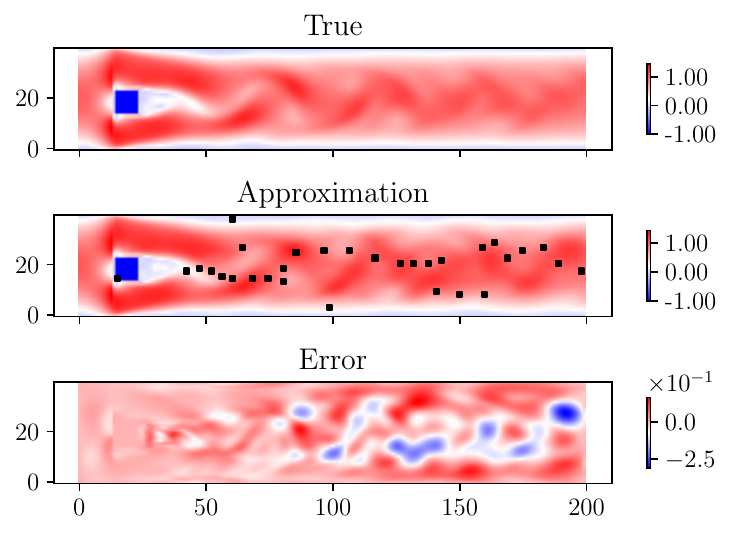}
		\caption[b]{Q-DEIM}
	\end{subfigure}
	\caption{t-Q-DEIM and Q-DEIM approximation of the thermal equation at $\p := 135$ at the time $t = 75.5$ seconds; each figure plots from top to bottom the true solution, the approximation, and the pointwise errors.}
	\label{fig:ns_tqdeim_qdeim_test_soln}
\end{figure}

The parameter set is obtained by selecting $N_{\text{s}} = 60$ equally-spaced samples from the parameter domain. This set is then further divided into the training and testing sets, with each containing $\ntrain = 45$ and $\ntest = 15$ randomly selected samples, respectively. Solving the Navier-Stokes equation at the sampled parameters, we obtain the training data $\cFtrain \in \R^{8000 \times 45 \times 320}$ and the testing data $\cFtest \in \R^{8000 \times 15 \times 320}$. Owing to the wide range of the Reynolds number (including vortex shedding behaviour), this represents a challenging problem for any sparse approximation technique.
\paragraph{Performance of t-Q-DEIM and Q-DEIM} Applying~\Cref{alg:tqdeim} to the training data $\cFtrain$, we obtain the t-Q-DEIM approximation. We also apply the Q-DEIM approach to the same dataset using~\Cref{alg:qdeim}. We first perform a sensitivity analysis to gauge the convergence of the approximation errors for either approach as a function of the reduced dimension $n$. This is plotted in~\Cref{fig:ns_sensitivty_analysis_abserr}. First, we observe that the magnitude of the error $\epsilon_{\text{abs}}$ is large. This is mainly an artefact of using the spectral norm and the large size of the problem. We repeat the same exercise but instead with the mean relative error defined as
$$
\epsilon_{\text{rel}} := \dfrac{1}{N_{\text{samples}}} \sum\limits_{i=1}^{N_{\text{samples}}} \dfrac{\| \mathfrak{f}(\p) - \widetilde{\mathfrak{f}}(\p) \|}{\| \mathfrak{f}(\p) \|}
$$
with $\p$ taken from the training set or the testing set and $N_{\text{samples}}$ either $N_{\text{train}}$ or $N_{\text{test}}$ as before.
\begin{figure}[h!]
	\centering
	\includegraphics[scale=0.6]{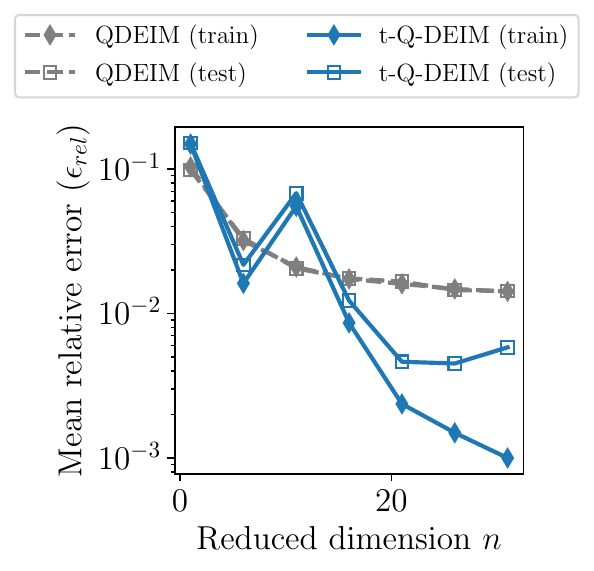}
	\caption{Navier-Stokes equation: Mean relative estimated errors over training (testing) set vs. $n$.}
	\label{fig:ns_sensitivty_analysis_relerr}
\end{figure} 
It can be observed that t-Q-DEIM approach yields a better approximation overall. Indeed, at two values of $n$, viz., $n = 1$ and $n = 11$, both the training and testing set relative errors for t-Q-DEIM are larger than those of Q-DEIM. Nevertheless, the overall trend reflects a superior performance of the proposed t-Q-DEIM approach. Next, we check the validity of the proposed error bound in~\Cref{thm:main_theorem}. For this, as done for previous examples, we plot the true error in the t-Q-DEIM approximation and the evaluated error bound in~\Cref{fig:ns_verify_error_bound}. Note that the bounds are plotted with the absolute t-spectral norm in the y-axis, and hence the large magnitude. Despite the larger t-spectral norm error, the upper bound holds, as shown for two choices of $n$, viz., $n = 21, 31$.
Finally, in~\Cref{fig:ns_tqdeim_qdeim_test_soln}, the true solution and the approximate solution with $n = 31$ at a test Reynolds number $Re = 135$, and time $t = 75.5$s are plotted. Both t-Q-DEIM and Q-DEIM approximations are compared. The pointwise errors support the conclusions in~\Cref{fig:ns_sensitivty_analysis_relerr}; the t-Q-DEIM errors are about an order of magnitude lower compared to the Q-DEIM pointwise errors. The figures showing the approximate solutions also show the locations of the sampling indices selected by the corresponding method. While the Q-DEIM indices display a more logical distribution, with sensors spread along the downstream vortices, some of the t-Q-DEIM sensor locations are characterised by a somewhat non-intuitive distribution. Up to four sensors are located, tightly clustered, at the very edge of the right boundary. 

This example again reflects the trend regarding the computational timing observed in the previous examples. The offline training cost for t-Q-DEIM is smaller by a very larger margin, compared to the respective cost for Q-DEIM. But, Q-DEIM is more beneficial with respect to the inference cost. 

\subsection{BMI experimental data}
\label{num:bmi}
\begin{figure}[t!]
	\centering
	\includegraphics[scale=0.65]{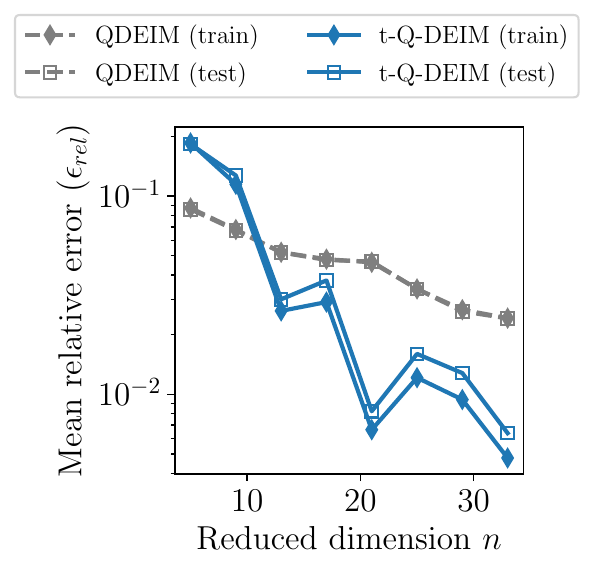}
	\caption{BMI example: Mean relative estimated errors over training (testing) set vs. $n$.}
	\label{fig:bmi_sensitivty_analysis}
\end{figure}

\begin{figure}[t!]
	\centering
	\begin{subfigure}[h!]{0.45\linewidth}
		\centering
		\includegraphics[scale=0.7]{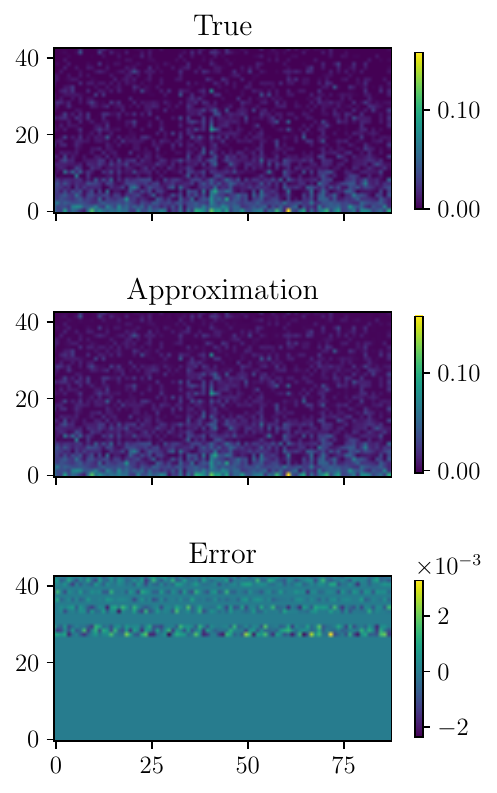}
		\caption[a]{t-Q-DEIM}
	\end{subfigure}\hfil
	\begin{subfigure}[h!]{0.45\linewidth}
		\centering
		\includegraphics[scale=0.7]{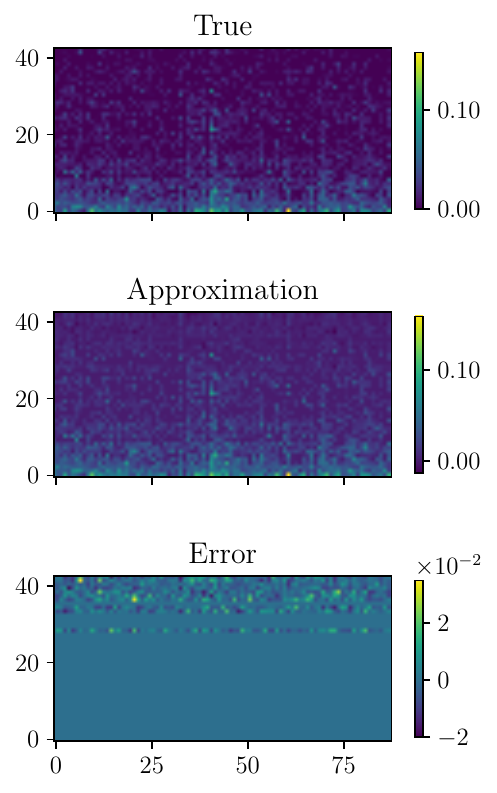}
		\caption[b]{Q-DEIM}
	\end{subfigure}
	\caption{t-Q-DEIM and Q-DEIM approximation of the BMI dataset at the time instance $t = 196$; each figure plots from top to bottom the true solution, the approximation, and the pointwise errors.}
	\label{fig:bmi_tqdeim_qdeim_test_soln}
\end{figure}

As our concluding example, we consider the tensor-valued dataset generated via a brain machine interface connected to the brain of a primate~\citep{Vyasetal18}. The dataset is taken from~\citep{tensor_repo}. The neuron activation responses of a primate are measured over time, in $88$ different trials. The activations are measured through sensors that monitor the activity of $43$ neurons. A  total of $200$ time steps are considered. The resulting third-order data tensor $\cF \in \R^{43 \times 200 \times 88}$ is then organized into training and test data tensors by randomly dividing the time samples. This yields the training data $\cFtrain \in \R^{43 \times 170 \times 88}$ and the test data $\cFtest \in \R^{43 \times 30 \times 88}$.
\paragraph{Performance of t-Q-DEIM and Q-DEIM}
We first carry out a sensitivity analysis to compare the approximation properties of t-Q-DEIM and Q-DEIM. We plot the relative mean error $\epsilon_{\text{rel}}$ as a function of the reduced dimension $n$; see~\Cref{fig:bmi_sensitivty_analysis}. The Q-DEIM approach seems to perform better at the initial two values of $n$, viz. $n = 5, 9$. However, after that, the performance of t-Q-DEIM is much better. This shows that the t-Q-DEIM approach we propose is able to deliver consistent performance even for datasets not originating from a PDE. Next, we show the performance of t-Q-DEIM and Q-DEIM on a test sample in~\Cref{fig:bmi_tqdeim_qdeim_test_soln}. It can be noted that, consistent with the previous examples, t-Q-DEIM yields better pointwise errors. Another key point to note here is that the BMI dataset consists of non-negative values. Therefore, a sparse approximation has to preserve this aspect. We note that both in the case of t-Q-DEIM and Q-DEIM some approximated values are negative. However, t-Q-DEIM results in far fewer negative values than Q-DEIM. Nevertheless, extending the t-Q-DEIM methodology to efficiently preserve the positivity of the data offers scope for future investigations.

\section{Conclusions and Outlook}
\label{sec:conclusion}
Accurate function reconstruction based on sparsely measured data is a problem with high practical relevance. When dealing with tensor-valued datasets (such as those arising from the numerical solution of parametrized PDEs), the approach of vectorizing/matricizing the dataset adopted by existing sparse approximation approaches such as the discrete empirical interpolation method leads to poorer approximation and/or predictions. We have presented an extension of the discrete empirical interpolation method which is able to handle third-order tensor-valued data, without the need for matricizing the data. Our proposed method, the t-Q-DEIM, avoids loss of geometric information in the tensor, resulting in improved approximations. We have developed rigorous, computable error bounds for the approximation error resulting from the t-Q-DEIM method. We further presented efficient computational strategies for obtaining the quantities needed in the construction of the t-Q-DEIM approximation, leveraging tools from the tensor t-product algebra. Additionally, the use of the t-SVD based on the FFT makes t-Q-DEIM a computationally attractive technique. Illustrations on the several numerical examples consisting of both, parametrized spatio-temporal data and heterogeneous 3-D data, show up to several orders of magnitude improvement in the approximation quality offered by the t-Q-DEIM over existing approaches. 

As promising future research directions, we highlight the need to investigate an extension of t-Q-DEIM that preserves certain properties of the dataset, such as its positivity. Furthermore, extensions to deal with noisy data will also make the t-Q-DEIM more relevant for use in practical applications.


\addcontentsline{toc}{section}{References}
\bibliographystyle{plainnat}
\bibliography{mybib}

\appendix
\section{Appendix}
\subsection{Computation of the t-pQR decomposition}
\label{appendix:tpqr}
The extension of the classical pivoted QR decomposition to the tensor t-product algebra was first proposed in~\citep{Haoetal13}. Consider the tensor $\cA \in \R^{m \times \ell \times q}$. Its tensor t-product QR decomposition is given by
\begin{align*}
	\cA * \cP = \cQ * \cR
\end{align*}
where $\cQ \in \R^{m \times m \times q}$ is an orthogonal tensor, $\cR \in \R^{m \times \ell \times q}$ is an f-upper triangular tensor (meaning each of its frontal slices is an upper triangular matrix), and $\cP \in\R^{\ell \times \ell \times q}$ is a permutation tensor. The computation of the t-pQR decomposition is as follows:

\begin{enumerate}
	\item Compute the Fourier domain representation of $\cA$:
	$$
		\chA := \texttt{fft}(\cA,\,[\,\,],\, 3) \in \R^{m \times \ell \times q}
	$$
	\item Compute the matrix pivoted QR decomposition of the first frontal slice of $\chA$:
	$$
		[\mathbf{Q}^{(1)},\,\mathbf{R}^{(1)}, \mathbf{P}^{(1)}] = \texttt{qr}(\chA(:,\,:,\,1))
	$$
	where $\mathbf{Q}^{(1)} \in \C^{m \times m}$ is an Hermitian matrix, $\mathbf{R}^{(1)} \in \C^{m \times \ell}$ is an upper triangular matrix, and the permutation matrix $\mathbf{P}^{(1)} \in \R^{\ell \times \ell}$ 
	\item For every remaining frontal slice of $\cA(:,\,:,\,j)$, $j \in \{2, 3,\ldots,q\}$, perform the matrix QR decomposition of its $\mathbf{P}^{(1)}$-column-permuted version, i.e.,
	$$
		[\mathbf{Q}^{(j)},\,\mathbf{R}^{(j)}, \sim] = \texttt{qr}\left(\chA(:,\,:,\,j) \mathbf{P}^{(1)} \right)
	$$
	\item Form tensors $\widehat{\cQ} \in \C^{m \times m \times q}, \widehat{\cR} \in \C^{m \times \ell \times q}$ and $\widehat{\cP} \in \R^{\ell \times \ell \times q}$ as follows:
		\begin{itemize}
			\item The $j$-th frontal slice of $\widehat{\cQ}$ consists of the matrix $\mathbf{Q}^{(j)}$,
			\item The $j$-th frontal slice of $\widehat{\cR}$ consists of the matrix $\mathbf{R}^{(j)}$, and
			\item Every frontal slice of $\widehat{\cP}$ consists of the same matrix $\mathbf{P}^{(1)}$.
		\end{itemize}
	\item Performing an inverse FFT on $\widehat{\cQ}, \widehat{\cR}$, and $\widehat{\cP}$ results in the tensors $\cQ, \cR$, and $\cP$:
	\begin{itemize}
		\item $\cQ := \texttt{ifft}\left(\widehat{\cQ},\,[\,\,],\,3\right)$,
		\item $\cR := \texttt{ifft}\left(\widehat{\cR},\,[\,\,],\,3\right)$, and
		\item $\cP := \texttt{ifft}\left(\widehat{\cP},\,[\,\,],\,3\right)$.
	\end{itemize}
\end{enumerate}
The resulting tensor $\cQ$ is orthogonal and $\cR$ is f-upper triangular.

For the sampling stage of the t-Q-DEIM algorithm~\Cref{alg:tqdeim}, we are only interested in the leading $n$ pivot indices of the first frontal slices, viz., the indices denoted by $\mathbf{P}^{(1)}$. Therefore, we terminate with Step 2 above.

\subsection{t-product of two orthogonal tensors}
\label{appendix:torth-svd}
We consider the orthogonal tensor $\cA \in \R^{\ell \times m \times n}$. Let $\cB \in \R^{\ell \times k \times n}$ be the tensor consisting of the first $k$ lateral slices of $\cA$. By definition, $\cB$ is an orthogonal tensor as well.

From the definition of the t-orthogonality~\Cref{def:t-orth}, it holds for the tensor t-product  $\cB^{\tpose} * \cA\left(:,\,1:k,\,:\right)$ that
$$
 \cB^{\tpose} * \cA\left(:,\,1:k,\,:\right) = \cI \in \R^{k \times k \times n}
$$.

Next, let us consider the tensor t-product $\cB^{\tpose} * \cA\left(:,\,k+1:m,\,:\right)$. For ease, let us define $\mathcal{Z} := \cB^{\tpose}$. The product can be expressed as
$$
\cB^{\tpose} * \cA\left(:,\,k+1:m,\,:\right) = \sum\limits_{i=1}^{\ell} \cZ\left(:,\,i,\,:\right) * \cA(i,\,:,\,:).
$$
Since the last $(m-k)$ lateral slices of $\cA$ are orthogonal to the first $k$ lateral slices of $\cB$, using \cref{eq:t-orth-property}, the above expression becomes
$$
\cB^{\tpose} * \cA\left(:,\,k+1:m,\,:\right) = \sum\limits_{i=1}^{\ell} \cZ\left(:,\,i,\,:\right) * \cA(i,\,:,\,:) = \mathcal{O} \in \R^{k \times (m-k) \times n}
$$
where $\mathcal{O}$ is a tensor of zeros. For the product $\cB^{\tpose} * \cA$, we then have
\begin{align*}
\cB^{\tpose} * \cA = \cB^{\tpose} * \left[\cB \qquad \cA\left(:,\,k+1 : m,\,:\right)\right] = \left[\mathcal{I}\qquad\mathcal{O}\right] \in \R^{k \times m \times n}
\end{align*}
In essence, this product results in a tensor, whose first $k$ lateral slices constitute an identity tensor and whose last $(m-k)$ lateral slices form a zero tensor.

\end{document}